\documentclass[10pt]{article}
\usepackage{amsmath,amssymb,amsthm}
\usepackage{color}
\usepackage[nohead,margin=1.0in]{geometry}
\usepackage{booktabs}
\usepackage{threeparttable}
\usepackage{algorithm,algorithmic}
\usepackage{verbatim}
\usepackage{url}
\usepackage{graphicx}
\usepackage{epstopdf}
\usepackage{mathtools}
\usepackage{subfig}
\graphicspath{{./pics/}}
\allowdisplaybreaks

\usepackage{array}
\theoremstyle{plain}
\newtheorem{theorem}{Theorem}[section]
\newtheorem{remark}{Remark}[section]

\newtheorem{lemma}{Lemma}[section]

\newtheorem{proposition}{Proposition}[section]
\newtheorem{corollary}{Corollary}[section]
\newtheorem{example}{Example}[section]

\numberwithin{equation}{section}

\renewcommand{\d}{\mathrm{d}}
\def\II{(\Omega)}

\title{Recovering the Potential and Order in One-Dimensional Time-Fractional Diffusion with Unknown Initial Condition and Source\thanks{The work of B. Jin is supported by UK EPSRC grant EP/T000864/1, and that of Z.Z. by Hong Kong RGC grant (No. 15304420).}
}
\author{Bangti Jin\thanks{Department of Computer Science, University College London, Gower Street, London WC1E 6BT, UK (\texttt{b.jin@ucl.ac.uk, bangti.jin@gmail.com})} \and Zhi Zhou\thanks{Department of Applied Mathematics, The Hong Kong Polytechnic University, Kowloon, Hong Kong. (\texttt{zhizhou@polyu.edu.hk})}}

\begin{document}

\maketitle
\begin{abstract}
This paper is concerned with an inverse problem of recovering a potential term and fractional order in
a one-dimensional subdiffusion problem, which involves a Djrbashian-Caputo fractional derivative of order
$\alpha\in(0,1)$ in time, from the lateral Cauchy data. In the model, we do not assume a full knowledge
of the initial data and the source term, since they might be unavailable in some practical applications.
We prove the unique recovery of the spatially-dependent potential coefficient and the order $\alpha$ of
the derivation simultaneously from the measured trace data at one end point, when the model is equipped
with a boundary excitation with a compact support away from $t=0$. One of the initial data and the source can also
be uniquely determined, provided that the other is known. The analysis employs a representation of the solution and the time analyticity of
the associated function. Further, we discuss a two-stage procedure, directly inspired by the analysis, for
the numerical identification of the order and potential coefficient, and illustrate the feasibility of the
recovery  with several numerical experiments.
\end{abstract}

\section{Introduction}
This paper is concerned with an inverse problem associated with a one-dimensional time-fractional diffusion
equation. Let $\Omega=(0,1)$ be the unit interval, and $T>0$ be a fixed final time.
Consider the following time-fractional diffusion problem for the function $u$:
\begin{equation}\label{eqn:fde}
  \left\{\begin{aligned}
     \partial_t^\alpha u - \mathcal{A}u  & = f, \quad\mbox{in }\Omega\times(0,T],\\
     -a(0)\partial_xu(0,t) & = g(t), \quad \mbox{in }(0,T],\\
     u(1,t) & = 0,\quad \mbox{in } (0,T],\\
     u(0) &= u_0, \quad \mbox{in }\Omega,
  \end{aligned}\right.
\end{equation}
where $f\in L^2(\Omega)$ and $u_0\in L^2(\Omega)$ are the (unknown) source term and initial
data, respectively. The second-order elliptic operator $\mathcal{A}$ is defined by
\begin{equation*}
  \mathcal{A}u(x) = \partial_x(a(x)\partial_x u(x)) - q(x)u(x),\quad x\in\overline{\Omega},
\end{equation*}
for $a\in C^2(\overline{\Omega})$ and $q\in L^\infty({\Omega})$, and $a\geq a_0>0$ and $q\geq 0$ in
$\Omega$. In the model, the notation $\partial_t^\alpha u$ denotes the Djrbashian-Caputo fractional
derivative of order $\alpha\in (0,1)$ in time, defined by \cite{Podlubny:1999,KilbasSrivastavaTrujillo:2006,Jin:2021}
\begin{equation*}
  \partial_t^\alpha u (t) =
    \frac{1}{\Gamma(1-\alpha)}\int_0^t(t-s)^{-\alpha}u'(s)\d s.
\end{equation*}
It is known that the Djrbashian-Caputo fractional derivative $\partial_t^\alpha u(t)$ recovers
the usual first-order derivative $u'(t)$ as $\alpha\to1^-$, and accordingly the model
\eqref{eqn:fde} reduces to the classical diffusion equation.

The model \eqref{eqn:fde} has been studied extensively in the engineering, physical and mathematical literature
due to its extraordinary capability for describing anomalous diffusion phenomena. It
can be viewed as the macroscopic counterpart of continuous time random walk in which the waiting time
between consecutive particle jumps follows a heavy-tailed distribution with a divergent mean, and the probability
density function of the particle appearing at location $x$ at time $t>0$ satisfies a model of the
form \eqref{eqn:fde}, in analogy with the classical diffusion equation for Brownian motion \cite{MetzlerKlafter:2000}. The model \eqref{eqn:fde}
inherits certain analytic properties of the latter, but also differs considerably
due to the presence of the nonlocal fractional derivative term $\partial_t^\alpha u$: it has limited smoothing
property in space and slow asymptotic decay in time \cite{SakamotoYamamoto:2011,Jin:2021}. The list of
successful applications is long and still fast growing, including thermal diffusion in fractal domains \cite{Nigmatulin:1986}, dispersion in a heterogeneous aquifer \cite{AdamsGelhar:1992} and
transport in column experiments \cite{HatanoHatano:1998} etc. See the comprehensive reviews \cite{MetzlerKlafter:2000,MetzlerJeon:2014} for the derivation of relevant mathematical models and many applications in physics and biology.

In this work, the inverse problem of interest is to recover the potential
$q\in\mathcal{Q}:=\{q\in C(\overline\Omega): q\geq 0\ \mbox{in }\Omega\}$ in the elliptic operator
$\mathcal{A}$ and the order $\alpha$ of derivation from the boundary observational data at the left end point
$h(t)=u(0,t)$ for $t\in[0,T]$. Note that in the model \eqref{eqn:fde}, besides the potential
$q$ and the order $\alpha$, the space-dependent source term $f$ and the initial data $u_0$ are both unknown.
The situation that the initial data $u_0$ is inaccessible arises naturally, e.g., in heat conduction
in high-temperature furnace \cite{WangChengNakagawaYamamoto:2010}. To make the matter worse, only one
single boundary observation data (at the left end) is available, which makes the inverse problem
much more challenging both mathematically and numerically. The ability of
choosing the boundary excitation $g(t)$ is indispensable for the unique recovery, without which the
identifiability generally does not hold, as indicated by example \ref{exam:counter-unique} below. In theorems
\ref{thm:init} and \ref{thm:main}, we present a uniqueness result for recovering the potential $q$ in
the operator $\mathcal{A}$ and the order $\alpha$ of derivation. The proof employs suitable
solution representation in proposition \ref{prop:solrep}, analyticity in time in proposition
\ref{prop:analyticity} and Gel'fand-Levitan
theory. Further, we discuss the numerical reconstruction by a two-stage procedure inspired directly
by the analysis: at stage i, we numerically continuate the boundary observation data $h(t)$ by rational
functions, and at stage ii, we perform a standard least-squares procedure to recover the potential
$q$ using the conjugate gradient method (with proper early stopping). The simulation study with exact
data indicates that the recovery is feasible. The uniqueness result, the two-stage recovery
procedure and the numerical verification of the recovery represent the main contributions of this work.

Next we situate the work in existing literature. The recovery of the space-dependent potential $q$
in the classical diffusion equation from lateral Cauchy data has been extensively discussed, and
several uniqueness results have been obtained \cite{Pierce:1979,Marayama:1981,SuzukiMurayama:1980}.
The study on related inverse problems for time-fractional models is of more recent origin, starting from
\cite{ChengNakagawa:2009} (see \cite{JinRundell:2015} for an early tutorial)  and there are a few works
on recovering a spatially dependent potential from lateral Cauchy data \cite{RundellYamamoto:2018,RundellYamamoto:2020,
WeiYan:2021,JingYamamoto:2021}. Rundell and Yamamoto \cite{RundellYamamoto:2018} showed
that the lateral Cauchy data can uniquely determine the spectral data when $u_0\equiv f\equiv
0$, and proved the uniqueness of the potential $q$ by the classical Gel'fand-Levitan theory. They also
proposed a recovery procedure based on Newton's method and empirically studied the singular value
spectrum of the linearized forward map, showing the severe ill-posed nature of the inverse problem. Later,
they \cite{RundellYamamoto:2020} relaxed the regularity condition on the boundary excitation
$g(t)$ (in a suitable Sobolev space in time). Recently, Jing and Yamamoto \cite{JingYamamoto:2021}
proved the identifiability of multiple parameters (including order, spatially dependent potential,
initial value and Robin coefficients in the boundary condition) simultaneously in the one-dimensional
subdiffusion / diffusion-wave (i.e., $\alpha\in(0,2)$) equation with a zero boundary condition and source,
excited by a nontrivial initial condition from the lateral Cauchy data at both end points (cf. remark
\ref{rmk:JingYamamoto} for details). See also the work \cite{WeiYan:2021}
for relevant results in the diffusion wave case; and \cite{JingPeng:2020} for the case of a Robin
boundary conditions. In all these existing works,
the initial condition / source is assumed to be fully known, so that the forward map is well defined, which
differs from the current work. There are two closely related inverse problems to the concerned one. (i)
is to recover the spatially dependent potential $q$ from the terminal
data $u(T)$ \cite{ZhangZhou:2017,KaltenbacherRundell:2019,JinZhou:2021ip},
which enjoys much better stability estimates (e.g., local Lipschitz stability) and effective
iterative algorithms for numerical recovery, e.g., fixed point iterations. {(ii) is to recover
a time-dependent potential $q(t)$ from time-dependent observations \cite{JiangWei:2021} (or a time-dependent
source from observation at one point \cite{LiuWenZhang:2020}), which behaves similarly to (i) due to the
directional alignment of the unknown and observations.}

The rest of the paper is organized as follows. In section \ref{sec:prelim}, we derive a crucial
representation of the solution to the direct problem \eqref{eqn:fde}. Then in section \ref{sec:unique}, we prove the unique recovery of
the order $\alpha$ and the potential $q$. In section \ref{sec:alg}, we describe a two-stage
numerical algorithm for recovering the potential $q$. Last, we present several numerical
experiments {to show} the feasibility of the simultaneous recovery in section \ref{sec:numer}.
Throughout, the notation $c$ denotes a generic constant which may differ at each occurrence, and
$(\cdot,\cdot)$ denotes the standard $L^2(\Omega)$ inner product (or duality pairing).

\section{Well-posedness of the direct problem}\label{sec:prelim}
In this section, we collect several preliminary results on the direct problem \eqref{eqn:fde},
especially the solution representation, which will play an important role in the study.

\subsection{Preliminaries}
First we describe several preliminary results that will be used in deriving the solution
representation. We use extensively the two-parameter Mittag-Leffler function
$E_{\alpha,\beta}(z)$ defined by (see, e.g., \cite{GorenfloKilbasMainardi:2020}
and \cite[Section 3.1]{Jin:2021})
\begin{equation*}
  E_{\alpha,\beta}(z) = \sum_{k=0}^\infty \frac{z^k}{\Gamma(k\alpha+\beta)},\quad z\in\mathbb{C}.
\end{equation*}
It is an entire function of order $\frac{1}{\alpha}$ and type one.
The following properties hold \cite[Section 3.1]{Jin:2021}.

\begin{lemma}\label{lem:mlf}
For any $\alpha\in (0,2)$ and $\beta\in\mathbb{R}$, the following statements hold.
\begin{itemize}
 \item[{\rm(i)}] For any $\varphi\in(\frac{\alpha}{2}\pi,\min(\pi,\alpha\pi))$,  the following asymptotics hold
\begin{equation*}
  |E_{\alpha,\beta}(z)|\leq \left\{\begin{aligned}
    c(1+|z|)^{-1},& \quad \beta-\alpha\not\in \mathbb{Z}_-\cup\{0\},\\
    c(1+|z|^2)^{-1},& \quad \beta-\alpha\in \mathbb{Z}_-\cup\{0\},
  \end{aligned}\right. \quad \varphi \leq |\arg(z)|\leq \pi.
\end{equation*}
\item[{\rm(ii)}] For any $\lambda>0$, the following Laplace transform relation holds
\begin{equation*}
  \mathcal{L}[t^{\beta-1}E_{\alpha,\beta}(-\lambda t^\alpha)](z) = \frac{z^{\beta-\alpha}}{\lambda + z^\alpha},\quad \Re(z)>0.
\end{equation*}
\item[{\rm(iii)}] The following differentiation formula holds for any $\lambda\in\mathbb{C}$
\begin{align*}
  \frac{\d}{\d t}E_{\alpha,1}(-\lambda t^\alpha) &= -\lambda t^{\alpha-1}E_{\alpha,\alpha}(-\lambda t^\alpha).
\end{align*}
\end{itemize}
\end{lemma}

Next, we introduce Bochner-Sobolev spaces $W^{\alpha,p}(0,T;X)$, for a UMD space $X$ (see \cite[Chapter 4]{Hytonen:2016}
for the definition of UMD spaces, which include Sobolev spaces $W^{s,p}(\Omega)$ with $s\ge0$ and $1<q<\infty$).
For any $s\ge 0$ and $1\le p< \infty$, we denote by $W^{s,p}(0,T;X)$ the space of functions $v:(0,T)\rightarrow X$,
with the norm defined by complex interpolation. Equivalently, the space is equipped with the quotient norm
\begin{align*}
\|v\|_{W^{s,p}(0,T;X)}&:= \inf_{\widetilde v}\|\widetilde v\|_{W^{s,p}({\mathbb R};X)}
:= \inf_{\widetilde v} \|  \mathcal{F}^{-1}[ (1+|\xi|^2)^{\frac{s}{2}} \mathcal{F} (\widetilde v)(\xi) ]\|_{L^p(\mathbb{R};X)},
\end{align*}
where the infimum is taken over all possible $\widetilde v$ that extend $v$ from $(0,T)$ to ${\mathbb R}$,
and $\mathcal{F}$ denotes the Fourier transform.
In case that $X=\mathbb{R}$, we denote $W^{s,p}(0,T;\mathbb{R})$ by $W^{s,p}(0,T)$ for convenience.
The next lemma provides a norm equivalence result \cite[Lemma 2.3]{JinKianZhou:2021}.
\begin{lemma}\label{lem:norm-equiv}
Let $\alpha\in(0,1)$ and $p\in [1,\infty)$ with $\alpha p>1$.
\begin{itemize}
\item[{\rm(i)}] If $v(0)=0$  and $ \partial_t^\alpha v \in L^p(0,T;X)$, then  $v \in {W^{\alpha,p}(0,T;X)}$ and
\begin{align*}
 \|v\|_{W^{\alpha,p}(0,T;X)} \le c \|  \partial_t^\alpha v \|_{L^p(0,T;X)}.
 \end{align*}
\item[{\rm(ii)}] If $v (0)=0$ and $ v \in  W^{\alpha,p}(0,T;X)$, then $\partial_t^\alpha v \in L^p(0,T; X)$
and
\begin{align*}
\|  \partial_t^\alpha v \|_{L^p(0,T;X)} \le c  \|v\|_{W^{\alpha,p}(0,T;X)}.
\end{align*}
\end{itemize}
\end{lemma}

\subsection{Well-posedness of the direct problem}

Now we study the direct problem \eqref{eqn:fde}, especially the solution representation. One distinct feature of
problem \eqref{eqn:fde} is that it involves a nonzero Neumann boundary condition, which has not been
extensively studied in the literature (\cite{RundellYamamoto:2020,KianYamamoto:2021} for relevant works). Following
\cite{RundellYamamoto:2020}, we exploit the one-dimensional nature of problem \eqref{eqn:fde}, and derive a series
representation of the solution $u$. The derivation is based on the standard separation of variable technique (see, e.g.,
\cite{SakamotoYamamoto:2011} and \cite[Section 6.2]{Jin:2021}). Specifically, let $A$ be the realization of the
elliptic operator $-\mathcal{A}$ in $L^2(\Omega)$, with its domain
$$D(A):=\{v\in L^2(\Omega): \mathcal{A}v\in L^2(\Omega), v'(0)=v(1)=0\}.$$
Let $\{(\lambda_n,\varphi_n)\}_{n=1}^\infty$ be the eigenpairs of the operator $A$, i.e.,
\begin{equation}\label{eqn:eigen}
  \left\{\begin{aligned}
     -\mathcal{A}\varphi_n = \lambda_n\varphi_n,\quad \mbox{in }\Omega,\\
     \varphi_n'(0) = 0, \quad \varphi_n(1) = 0.
  \end{aligned}\right.
\end{equation}
By the standard Sturm-Liouville theory \cite{LevitanSargsjan:1975},
the spectrum of the operator $A$ consists of a strictly increasing sequence
of positive eigenvalues $\{\lambda_n\}_{n=1}^\infty$ and the associated
eigenfunctions $\{\varphi_n\}_{n=1}^\infty$ can be chosen to form an orthonormal basis
of the space $L^2(\Omega)$. By means of Liouville transformation, we deduce that the eigenvalues
$\lambda_n$ grow asymptotically as $O(n^2)$ \cite{LevitanSargsjan:1975} (also known from Weyl's law \cite{Weyl:1912}):
\begin{equation*}
  \lambda_n = \Big(\int_0^1a(x)^{-\frac12}\d x\Big)^{-1}(n-\tfrac12)^2\pi^2+O(1),\quad n\to \infty.
\end{equation*}
The unnormalized eigenfunctions $\varphi_n$
satisfy the following asymptotics \cite[Section 2 of Chapter 1]{LevitanSargsjan:1975}:
\begin{equation}\label{eqn:eigf-b}
\varphi_n(x) = \Big(\sqrt{2}\cos\Big[\frac{(n-\frac12)\pi }{\int_0^1 a(s)^{-\frac{1}{2}}\d s}\int_0^x a(s)^{-\frac12}\d s\Big] + O\Big(\frac1n\Big)\Big) a(x)^{-\frac{1}{4}},
\qquad n\in \mathbb{N}.
\end{equation}
This estimate implies that the $L^2(\Omega)$-orthonormal eigenfunctions are uniformly bounded.

Then we define the fractional power $A^s$, $s\geq0$, by
\begin{equation*}
  A^s v = \sum_{n=1}^\infty \lambda_n^s(v,\varphi_n)\varphi_n,
\end{equation*}
with its domain $\{v\in L^2(\Omega): A^sv\in L^2(\Omega)\}$, and the associated graph norm $\|\cdot\|_{D(A^s)}$ given by
\begin{equation*}
  \|v\|_{D(A^s)} = \Big(\sum_{n=1}^\infty \lambda_n^{2s}(v,\varphi_n)^2\Big)^\frac{1}{2}.
\end{equation*}
For $s<0$, $D(A^{s})$ is identified with the dual space of $D(A^{-s})$.

With these preliminaries, we can now study the direct problem \eqref{eqn:fde}. In view of the
linearity of the problem, we may split the solution $u$ into two parts:
$u=u_i+u_b$, with $u_i$ and $u_b$ solving respectively
\begin{equation*}
  \left\{\begin{aligned}
     \partial_t^\alpha u_i - \mathcal{A}u_i  & = f, \quad\mbox{in }\Omega\times(0,T],\\
     -a(0)\partial_xu_i(0,t) & = 0, \quad \mbox{in }(0,T],\\
     u_i(1,t) & = 0,\quad \mbox{in } (0,T],\\
     u_i(0) &= u_0, \quad \mbox{in }\Omega,
  \end{aligned}\right.
  \quad\mbox{and}\quad
  \left\{\begin{aligned}
     \partial_t^\alpha u_b - \mathcal{A}u_b  & = 0, \quad\mbox{in }\Omega\times(0,T],\\
     -a(0)\partial_xu_b(0,t) & = g, \quad\mbox{in }(0,T],\\
     u_b(1,t) & = 0,\quad \mbox{in } (0,T],\\
     u_b(0) &= 0, \quad \mbox{in }\Omega.
  \end{aligned}\right.
\end{equation*}

With the eigenexpansion $\{(\lambda_n,\varphi_n)\}_{n=1}^\infty$, the solution $u_i$
can be represented by (see e.g., \cite{SakamotoYamamoto:2011} and \cite[Section 6.2]{Jin:2021})
\begin{align}
  u_i(x,t) & = \sum_{n=1}^\infty (u_0,\varphi_n)E_{\alpha,1}(-\lambda_nt^\alpha)\varphi_n(x) + \sum_{n=1}^\infty
  \int_0^t(t-s)^{\alpha-1}E_{\alpha,\alpha}(-\lambda_n(t-s)^\alpha)(f,\varphi_n)\d s\varphi_n(x).\label{eqn:sol-ui-0}
\end{align}
When the source $f$ is time independent, by Lemma \ref{lem:mlf}(iii), we have
\begin{align}
  u_i(x,t)
           & = \sum_{n=1}^\infty (u_0,\varphi_n)E_{\alpha,1}(-\lambda_nt^\alpha)\varphi_n(x) + \sum_{n=1}^\infty
  \lambda_n^{-1}[1-E_{\alpha,1}(-\lambda_nt^\alpha)](f,\varphi_n)\varphi_n(x)\nonumber\\
           & = \sum_{n=1}^\infty \big([(u_0,\varphi_n) - \lambda_n^{-1}(f,\varphi_n)]E_{\alpha,1}(-\lambda_nt^\alpha)+\lambda_n^{-1}(f,\varphi_n)\big)\varphi_n(x).\label{eqn:sol-ui}
\end{align}

Further, we have the following \textit{a priori} estimate on the solution $u_i$.
\begin{lemma}\label{lem:reg}
If $u_0 \in D(A^\beta)$ with $\beta\in(0,1)$
and $f\in L^p(0,T;L^2(\Omega))$ with $\alpha p\in(1,(1-\beta)^{-1})$,
then there is a unique solution
$u_i \in W^{\alpha,p}(0,T;L^2(\Omega))\cap L^p(0,T;H^2(\Omega))$ such that
\begin{equation}\label{eqn:reg}
 \|  u_i \|_{W^{\alpha,p}(0,T;L^2(\Omega))} +  \|  u_i \|_{L^p(0,T;H^2(\Omega))} \le c \left(\| A^\beta u_0 \|_{L^2(\Omega)} + \| f \|_{L^p(0,T;L^2(\Omega))} \right).
\end{equation}
\end{lemma}
\begin{proof}
When $f\equiv0$ and $u_0\in D(A^\beta)$, there holds (see e.g. \cite[Section 6.2]{Jin:2021})
$$ \| \partial_t^\alpha u_i  \|_{L^2\II}  +  \|  \mathcal{A} u_i  \|_{L^2\II} \le c t^{-\alpha(1-\beta)} \| A^\beta u_0 \|_{L^2(\Omega)},   $$
which further implies for $p\in(0,(\alpha(1-\beta))^{-1})$
$$ \| \partial_t^\alpha u_i  \|_{L^p(0,T;L^2\II)}  +  \|  \mathcal{A} u_i  \|_{L^p(0,T;L^2\II)} \le c_T \| A^\beta u_0 \|_{L^2(\Omega)}.   $$
Then the elliptic regularity theory implies $u_i \in  L^p(0,T;H^2(\Omega))$. Meanwhile, lemma \ref{lem:norm-equiv} (i) leads to
$u_i-u_0 \in W^{\alpha,p}(0,T;L^2(\Omega))$, so is $u_i$.
For the case that $u\equiv0$ and $f\in L^p(0,T;L^2\II)$, we refer to \cite[Theorem 6.11]{Jin:2021} for a detailed proof.
\end{proof}

Next we turn to the a representation of the solution $u_b$. We need the following identity.
\begin{lemma}\label{lem:mlf-conv}
For $\alpha\in(0,1)$ and $g\in C^1[0,T]$,  there holds
\begin{align*}
  &\int_0^t(t-s)^{\alpha-1}E_{\alpha,\alpha}(-\lambda_n(t-s)^\alpha)\partial_s^\alpha g(s)\d s \\
 =&  g(t) -  g (0)E_{\alpha,1}(-\lambda_nt^\alpha)- \int_0^t\lambda_n(t-s)^{\alpha-1}E_{\alpha,\alpha}(-\lambda_n(t-s)^\alpha)g(s)\d s.
\end{align*}
\end{lemma}
\begin{proof}
A variant of the identity  can be found in \cite[Section 2]{RundellYamamoto:2020},
and we recap the proof only for completeness. We denote the integral on the
left hand side by $S$. Then changing the order of integration gives
\begin{align*}
  S &= \frac{1}{\Gamma(1-\alpha)}\int_0^t(t-s)^{\alpha-1}E_{\alpha,\alpha}(-\lambda_n(t-s)^\alpha)\Big(\int_0^s(s-\xi)^{-\alpha}g'(\xi)\d \xi\Big)\d s\\
   & = \frac{1}{\Gamma(1-\alpha)}\int_0^t g'(\xi)\Big(\int_\xi^t(t-s)^{\alpha-1}E_{\alpha,\alpha}(-\lambda_n(t-s)^\alpha)(s-\xi)^{-\alpha}\d s\Big)\d \xi\\
   & = \frac{1}{\Gamma(1-\alpha)}\int_0^t g'(\xi)\Big(\int_0^{t-\xi}\eta^{\alpha-1}E_{\alpha,\alpha}(-\lambda_n\eta^\alpha)(t-\xi-\eta)^{-\alpha}\d \eta\Big)\d \xi,
\end{align*}
where the last line follows from the change of variables $\eta=t-s$. Moreover,
by the definition of $E_{\alpha,\alpha}(-\lambda_nt^\alpha)$ and
applying termwise integration \cite[(3.5)]{Jin:2021}, we have
$\frac{1}{\Gamma(1-\alpha)}\int_0^t \eta^{\alpha-1}E_{\alpha,\alpha}(-\lambda_n\eta^\alpha)(t-\eta)^{-\alpha}\d \eta = E_{\alpha,1}(-\lambda_nt^\alpha)$.
Thus, by integration by parts and lemma \ref{lem:mlf}(iii), we obtain
\begin{align*}
S &= \int_0^t g'(\xi) E_{\alpha,1}(-\lambda_n(t-\xi)^\alpha)\d \xi\\
  &= [g(\xi) E_{\alpha,1}(-\lambda_n(t-\xi)^\alpha)]_{\xi=0}^{\xi=t} - \int_0^tg(\xi)\frac{\d}{\d\xi}E_{\alpha,1}(-\lambda_n(t-\xi)^\alpha)\d \xi\\
   &=  g(t) -  g (0)E_{\alpha,1}(-\lambda_nt^\alpha)- \int_0^t\lambda_n(t-s)^{\alpha-1}E_{\alpha,\alpha}(-\lambda_n(t-s)^\alpha)g(s)\d s.
\end{align*}
This completes the proof of the lemma.
\end{proof}

Now we can derive a representation of the solution $u_b$, corresponding to nonzero
boundary conditions. The derivation essentially exploits the one-dimensional nature
of the problem.
\begin{proposition}\label{prop:solrep}
Let $\alpha\in(0,1)$ and $p>4/(3\alpha)$ .
Suppose $g_1\in L^p(0,T)$, $g_2\in W^{\alpha,p}(0,T)$ and $g_2(0)=0$. Then the solution $u$ to the following initial boundary value problem
\begin{equation}\label{eqn:solrep-0}
  \left\{\begin{aligned}
     \partial_t^\alpha u - \mathcal{A}u  & = 0, \quad\mbox{in }\Omega\times(0,T],\\
     -a(0)\partial_xu(0,t) & = g_1(t), \quad \mbox{in }(0,T],\\
     u(1,t) & = g_2(t),\quad \mbox{in } (0,T],\\
     u(0) &= 0, \quad \mbox{in }\Omega,
  \end{aligned}\right.
\end{equation}
can be represented by
\begin{equation}\label{eqn:solrep}
\begin{aligned}
u(x,t) &=\sum_{n=1}^\infty\int_0^t(t-s)^{\alpha-1}E_{\alpha,\alpha}(-\lambda_n(t-s)^\alpha)g_1(s)\d s \varphi_n(0)\varphi_n(x)\\
    &\quad -\sum_{n=1}^\infty \int_0^t(t-s)^{\alpha-1}E_{\alpha,\alpha}(-\lambda_n(t-s)^\alpha)g_2(s)\d sa(1)\varphi_n'(1)\varphi_n(x).
  \end{aligned}
\end{equation}
\end{proposition}
\begin{proof}
The derivation proceeds by homogenizing the boundary conditions.
First, we assume $g_1,g_2\in C^1[0,T]$, which will be relaxed below. Let
\begin{equation*}
  v(x,t)=u(x,t)-\phi_0(x)g_1(t)-g_2(t),\quad \mbox{with }\phi_0(x)=(2a(0))^{-1}(x-1)^2.
\end{equation*}
Clearly, $\phi'_0(0)=-a(0)^{-1}$. Then direct computation shows
\begin{align*}
  \partial_t^\alpha v & = \partial_t^\alpha u(x,t) - \phi_0(x)\partial_t^\alpha g_1(t) - \partial_t^\alpha g_2(t) = \mathcal{A} v(x,t) + \tilde f(x,t),
\end{align*}
with
\begin{equation*}
 \tilde f(x,t) = \mathcal{A}(\phi_0(x)g_1(t) + g_2(t)) - \phi_0(x)\partial_t^\alpha g_1(t) - \partial_t^\alpha g_2(t).
\end{equation*}
Next we evaluate the boundary conditions:
\begin{align*}
  -a(0)v'(0,t) &= -a(0)(u'(0,t) - \phi_0'(0)g_1(t)) \\
  &= -a(0)u'(0,t)+ a(0)\phi'(0)g_1(t) = g_1(t) - g_1(t)=0,\\
  v(1,t) & =  u(1,t) - g_2(t) = 0.
\end{align*}
Thus, the function $v$ satisfies
\begin{equation*}
  \left\{\begin{aligned}
     \partial_t^\alpha v - \mathcal{A}v  & = \tilde f, \quad\mbox{in }\Omega\times(0,T),\\
     -a(0)\partial_xv(0,t) & = 0, \quad \mbox{in }(0,T],\\
     v(1,t) & = 0,\quad \mbox{in } (0,T],\\
     v(0) &= {-\phi_0g_1(0)}-g_2(0), \quad \mbox{in }\Omega.
  \end{aligned}\right.
\end{equation*}
By the construction, we have $v(0) \in D(A^s)$ with $s\in(0,\frac14)$ and $\tilde f \in L^\infty(0,T; L^2\II)$.
Since the function $v$ satisfies homogeneous boundary conditions, by the standard separation of variable technique
(see e.g., \cite{SakamotoYamamoto:2011} or \cite[Section 6.1]{Jin:2021}), it can be represented by
\begin{align*}
  v(x,t) =& \sum_{n=1}^\infty
   E_{\alpha,1}(-\lambda_nt^\alpha)(\varphi_n,v (0))\varphi_n(x)\\
   & + \sum_{n=1}^\infty  \int_0^t(t-s)^{\alpha-1}E_{\alpha,\alpha}(-\lambda_n(t-s)^\alpha)(\varphi_n,\tilde f(s))\,\d s \,\varphi_n(x).
\end{align*}
Next we simplify the second term. By the definition of $\tilde f$ and using lemma \ref{lem:mlf-conv}, we deduce
\begin{align*}
  & \int_0^t(t-s)^{\alpha-1}E_{\alpha,\alpha}(-\lambda_n(t-s)^\alpha)(\varphi_n,\tilde f(s))\, \d s\\
  = & - \int_0^t(t-s)^{\alpha-1}E_{\alpha,\alpha}(-\lambda_n(t-s)^\alpha)[\partial_s^\alpha g_1(s)(\phi_0,\varphi_n)-g_1(s)(\mathcal{A}\phi_0,\varphi_n)]\, \d s\\
  & - \int_0^t(t-s)^{\alpha-1}E_{\alpha,\alpha}(-\lambda_n(t-s)^\alpha)[\partial_s^\alpha g_2(s)(1,\varphi_n)-g_2(s) (\mathcal{A}1,\varphi_n)\, \d s]\\
  = & - \Big[g_1(t)- g_1(0) E_{\alpha,1}(-\lambda_nt^\alpha)-\int_0^t\lambda_n(t-s)^{\alpha-1}E_{\alpha,\alpha}(-\lambda_n(t-s)^\alpha)g_1(s)\d s\Big](\phi_0,\varphi_n)\\
  & +\int_0^t(t-s)^{\alpha-1}E_{\alpha,\alpha}(-\lambda_n(t-s)^{\alpha})g_1(s)\, \d s (\mathcal{A}\phi_0,\varphi_n)\\
  & - \Big[g_2(t)- g_2(0) E_{\alpha, 1}(-\lambda_nt^\alpha)-\int_0^t\lambda_n(t-s)^{\alpha-1}E_{\alpha,\alpha}(-\lambda_n(t-s)^\alpha)g_2(s)\d s\Big](1,\varphi_n)\\
  &+ \int_0^t(t-s)^{\alpha-1}E_{\alpha,\alpha}(-\lambda_n(t-s)^{\alpha})g_2(s)\, \d s \, (\mathcal{A}1,\varphi_n)\\
  = & - \Big[g_1(t)- g_1(0) E_{\alpha, 1}(-\lambda_nt^\alpha)\Big](\phi_0,\varphi_n)
  - \Big[g_2(t)- g_2(0) E_{\alpha, 1}(-\lambda_nt^\alpha)\Big](1,\varphi_n)\\
  &+\int_0^t(t-s)^{\alpha-1}E_{\alpha,\alpha}(-\lambda_n(t-s)^\alpha)g_1(s)\, \d s
  \, \Big[\lambda_n(\phi_0,\varphi_n) + (\mathcal{A}\phi_0,\varphi_n)\Big]\\
  & +\int_0^t(t-s)^{\alpha-1}E_{\alpha,\alpha}(-\lambda_n(t-s)^\alpha)g_2(s)\, \d s
  \, \Big[\lambda_n(1,\varphi_n)+ (\mathcal{A}1,\varphi_n)\Big].
\end{align*}
Next we evaluate the last two terms in the square bracket. By Green's identity, for any $\phi\in H^2(\Omega)$,
\begin{align*}
  &\lambda_n(\phi,\varphi_n) + (\mathcal{A}\phi,\varphi_n)
  =-(\phi,\mathcal{A}\varphi_n) + (\mathcal{A}\phi,\varphi_n)
  =  a\phi'\varphi_n|_{x=0}^{x=1} - a\phi\varphi_{n}'|_{x=0}^{x=1}.
\end{align*}
Using this identity, the two terms in the square bracket can be evaluated as
\begin{align*}
  \lambda_n(\phi_0,\varphi_n) + (\mathcal{A}\phi_0,\varphi_n) & = -a(0)(-a(0)^{-1})\varphi_n(0) + a(0)(2a(0))^{-1}\varphi_{n}'(0) = \varphi_n(0),\\
  \lambda_n(1,\varphi_n) + (\mathcal{A}1,\varphi_n) & = -a\varphi_{n}'|_{x=0}^{x=1} = -a(1)\varphi_{n}'(1).
\end{align*}
Consequently, we arrive at
\begin{align*}
  & \int_0^t(t-s)^{\alpha-1}E_{\alpha,\alpha}(-\lambda_n(t-s)^\alpha)(\varphi_n,\tilde f(s))\d s\\
  = & - \Big[g_1(t)- g_1(0) E_{\alpha, 1}(-\lambda_nt^\alpha)\Big](\phi_0,\varphi_n)
  - \Big[g_2(t)- g_2(0) E_{\alpha, 1}(-\lambda_nt^\alpha)\Big](1,\varphi_n)\\
  &+\int_0^t(t-s)^{\alpha-1}E_{\alpha,\alpha}(-\lambda_n(t-s)^\alpha)g_1(s)\d s \varphi_n(0)\\
  & -\int_0^t(t-s)^{\alpha-1}E_{\alpha,\alpha}(-\lambda_n(t-s)^\alpha)g_2(s)\d sa(1)\varphi_n'(1).
\end{align*}
This identity and the representation of $v$ lead to
\begin{align*}
  v(x,t) =&-\phi_0(x)g_1(t) - g_2(t) +
  \sum_{n=1}^\infty\int_0^t(t-s)^{\alpha-1}E_{\alpha,\alpha}(-\lambda_n(t-s)^\alpha)g_1(s)\d s \varphi_n(0)\varphi_n(x)\\
  & -\sum_{n=1}^\infty \int_0^t(t-s)^{\alpha-1}E_{\alpha,\alpha}(-\lambda_n(t-s)^\alpha)g_2(s)\d sa(1)\varphi_n'(1)\varphi_n(x).
\end{align*}
Now the definition of the function $u$ concludes the proof of the proposition for $g_1,g_2\in C^1[0,T]$.

Next, we prove the assertion for $g_1 \in L^p(0,T),g_2\in W^{\alpha,p}(0,T)$ by a density argument.
First, we assume that  $g_1 \equiv 0$ and $g_2 \in W^{\alpha,p}(0,T) $.
Since $C^1[0,T]$ is dense in  $W^{\alpha,p}(0,T)$,
for any $ g_2\in W^{\alpha,p}(0,T)$, there exists a sequence $ \{g_2^\ell\}_{\ell=1}^\infty
\subset C^1[0,T]$ such that  $g_2^\ell (0)=g_2(0)$, and  $g_2^\ell\to g_2$
in $W^{\alpha,p}(0,T)$.
Then by lemma \ref{lem:norm-equiv} (ii) we conclude  $\partial_t^\alpha g_2^\ell\to \partial_t^\alpha g_2$
 in $L^2(0,T)$.
 Then by the \textit{a priori} estimate in lemma \ref{lem:reg}, we have
 $v_{g_2^\ell}\to v_{g_2}$ in $W^{\alpha,p}(0,T;L^2(\Omega))\cap L^p(0,T;H^2(\Omega))$.
Then,  we assume that  $g_2 \equiv 0$ and $g_1 \neq 0$. By the expression of the function $u$,
for any $g_1\in C^1[0,T]$, there holds for arbitrarily small $\epsilon>0$ \cite[Theorem 5.2]{JinKianZhou:2021}
\begin{equation*}
\| u \|_{W^{\alpha,p}(0,T; D(A^{-\frac14-\epsilon}))} + \| u \|_{L^p(0,T; D(A^{\frac34-\epsilon}))} \le c\| g_1 \|_{L^p(0,T)}.
\end{equation*}
By means of the complex interpolation and the density of $C^1[0,T]$ in  $L^p(0,T)$,
we deduce that for any $g_1\in L^p(0,T)$, the representation \eqref{eqn:solrep} provides a solution of problem
\eqref{eqn:solrep-0}, and $u \in W^{s,p}(0,T;L^2(\Omega))$ with $s<\frac{3\alpha}{4}$, which
embeds compactly into  $C([0,T];L^2(\Omega))$ \cite[Theorem 5.2]{Amann:2000}.
\end{proof}

The following representation of the observation data $h(t)=u(0,t)$ is direct.
\begin{corollary}\label{cor:sol-rep}
The observation data $h(t)=u(0,t)$ to problem \eqref{eqn:fde} can be represented by
\begin{equation*}
  h(t) = \sum_{n=0}^\infty \rho_nE_{\alpha,1}(-\lambda_nt^\alpha)  + \sum_{n=1}^\infty \int_0^t(t-s)^{\alpha-1}E_{\alpha,\alpha}(-\lambda_n(t-s)^\alpha)\varphi_n(0)g(s)\d s\varphi_n(x),
\end{equation*}
with
\begin{equation*}
  \lambda_0=0, \quad \rho_0 = \sum_{n=1}^\infty\lambda_n^{-1}(f,\varphi_n)\varphi_n(0),\quad  \rho_n(x) = [(u_0,\varphi_n)-\lambda_n^{-1}(f,\varphi_n)]\varphi_n(0), \quad n=1,2,\ldots
\end{equation*}
\end{corollary}
\begin{proof}
The assertion is direct from proposition \ref{prop:solrep} and the identity \eqref{eqn:sol-ui}. Indeed, we have
\begin{align*}
  u(x,t) &= \sum_{n=1}^\infty E_{\alpha,1}(-\lambda_nt^\alpha)(\varphi_n,u_0)\varphi_n(x) \\
     &\quad + \sum_{n=1}^\infty \int_0^t(t-s)^{\alpha-1}E_{\alpha,\alpha}(-\lambda_n(t-s)^\alpha)((f,\varphi_n)+\varphi_n(0)g(s))\d s\varphi_n(x).
\end{align*}
Then by lemma \ref{lem:mlf}(iii), we have $\int_0^ts^{\alpha-1}E_{\alpha,\alpha}(-\lambda_ns^\alpha)\d s = \lambda_n^{-1}(1- {E_{\alpha,1}}(-\lambda_nt^\alpha)).$ This directly leads to the desired identity.
\end{proof}

\begin{remark}
Let the solution operators $E(t)$ and $F(t)$ be defined by
\begin{equation*}
  E(t)v = \sum_{n=1}^\infty (v,\varphi_n)t^{\alpha-1}E_{\alpha,\alpha}(-\lambda_nt^\alpha)\varphi_n \quad
  \mbox{and}\quad F(t)v = \sum_{n=1}^\infty (v,\varphi_n)E_{\alpha,1}(-\lambda_nt^\alpha)\varphi_n.
\end{equation*}
Then the solution $u$ can be formally represented with
\begin{equation*}
  u(t) = F(t)u_0 + \int_0^t E(t-s)[\delta_0(x)g(s) + f] \d s,
\end{equation*}
where $\delta_0(x)$ denotes the Dirac delta function concentrated at
$x=0$. Indeed, we can expand $\delta_0(x)$ in terms of the $L^2(\Omega)$ orthonormal basis
$\{\varphi_n\}_{n=1}^\infty$ in $D(A)'$ {\rm(}the dual space of $D(A)${\rm)} as
$\delta(x) = \sum_{n=1}^\infty \varphi_n(0)\varphi_n(x)$.
Then substituting the identity and collecting the terms lead to the desired identity.
Alternatively, the solution can also be represented concisely using the fractional
$\theta$ functions; see \cite[Section 7.1]{Jin:2021} for relevant discussions.
\end{remark}

Next we study the convergence of the series in corollary \ref{cor:sol-rep}. First, note that for
$f\in D(A^{-s}) $ with $0\leq s<\frac{3}{4}$, the constant $\rho_0$ is well defined. Indeed, by \eqref{eqn:eigf-b}
and the Cauchy-Schwarz inequality, and asymptotics for the eigenvalues $\{\lambda_n\}_{n=1}^\infty$, we have
\begin{equation*}
  |\rho_0| \leq c\Big(\sum_{n=1}^\infty \lambda_n^{-2s}(f,\varphi_n)^2\Big)^\frac{1}{2}\Big(\sum_{n=1}^\infty \lambda_n^{-2(1-s)}\Big)^\frac12
  \leq c\|f\|_{D(A^{-s})} \Big(\sum_{n=1}^\infty n^{-4(1-s)}\Big)^\frac12 <\infty.
\end{equation*}
For the remaining series, we give two analyticity results concerning the following two auxiliary functions:
\begin{equation*}
  h_1(t) = \sum_{n=1}^\infty  \rho_n E_{\alpha,1}(-\lambda_nt^\alpha)\quad\mbox{and}\quad
h_2(t)=\sum_{n=1}^\infty |\varphi_n(0)|^2t^{\alpha-1}E_{\alpha,\alpha}(-\lambda_nt^\alpha).
\end{equation*}
They arise naturally in the uniqueness proof, and the analyticity plays an important role
in section \ref{sec:unique}.
\begin{proposition}\label{prop:analyticity}
For $u_0,f\in L^2(\Omega)$, the following statements hold.
\begin{itemize}
  \item[{\rm(i)}] Both $h_1(t)$ and $h_2(t)$ are analytic in $t$ on $(0,\infty)$.
  \item[{\rm(ii)}] The Laplace transforms of $h_1(t)$ and $h_2(t)$ exist and are given respectively by
  \begin{align*}
    \mathcal{L}[h_1(t)](z) = \sum_{n=1}^\infty \frac{\rho_nz^{\alpha-1}}{z^\alpha+\lambda_n}
    \quad\mbox{and}\quad \mathcal{L}[h_2(t)](z) = \sum_{n=1}^\infty \frac{\varphi_n(0)^2}{z^\alpha+\lambda_n},\quad \forall \Re(z)>0.
  \end{align*}
\end{itemize}
\end{proposition}
\begin{proof}
By lemma \ref{lem:mlf}(i), there exists a constant $c$ and $\theta_0>0$ such that
\begin{equation*}
  |E_{\alpha,1}(-\lambda_nz^\alpha)|\leq c\lambda_n^{-1} |z|^{-\alpha},\quad \forall n\in\mathbb{N}\mbox{ and } z\in \Sigma:=\{z\in \mathbb{C}: |\arg(z)|<\theta_0\}.
\end{equation*}
By the asymptotic estimate  \eqref{eqn:eigf-b} of the eigenfunctions $\{\varphi_n\}_{n=1}^\infty$, we deduce $|\varphi_n(0)| \le c$  uniformly in $n$. Consequently,
\begin{align*}
  |E_{\alpha,1}(-\lambda_nz^\alpha) \rho_n| \leq  c|z|^{-\alpha}\lambda_n^{-1}(|(u_0,\varphi_n)|+\lambda_n^{-1}|(f,\varphi_n)|),\quad \forall x\in \overline\Omega, z\in\Sigma,
\end{align*}
and by the Weyl's asymptotics of the eigenvalues $\lambda_n$ and the Cauchy-Schwarz inequality, we have
\begin{align}
 | h_1(z)|&\leq\sum_{n=1}^\infty |E_{\alpha,1}(-\lambda_nt^\alpha) \rho_n| \leq
   c|z|^{-\alpha}\sum_{n=1}^\infty \lambda_n^{-1}(|(u_0,\varphi_n)|+\lambda_n^{-1}|(f,\varphi_n)|)\nonumber\\
  & \leq c|z|^{-\alpha}\Big(\sum_{n=1}^\infty|(u_0,\varphi_n)|^2+\lambda_n^{-2}|(f,\varphi_n)|^2\Big)^\frac12\Big(\sum_{n=1}^\infty n^{-4}\Big)^\frac12<c|z|^{-\alpha},\quad \forall z\in \Sigma.\label{eqn:bdd-sum}
\end{align}
Since $E_{\alpha,1}(-\lambda_nz^\alpha)$ is analytic in $z\in \Sigma$, we deduce that the
series in analytic in $t$. The analyticity of $h_2(t)$
follows similarly as $h_1(t)$, but with the following estimate from lemma \ref{lem:mlf}(i):
there exists a constant $c$ and $\theta_0>0$ such that
\begin{equation*}
  |E_{\alpha,\alpha}(-\lambda_nz^\alpha)|\leq c(1+\lambda_n^2 |z|^{2\alpha})^{-1},\quad \forall n\in\mathbb{N}\mbox{ and } z\in \Sigma:=\{z\in \mathbb{C}: |\arg(z)|<\theta_0\}.
\end{equation*}
Then by Weyl's asymptotics of the eigenvalues $\lambda_n$, for any $\gamma\in (\frac12,1)$,
\begin{align*}
   |h_2(z)| & \leq c|z|^{\alpha-1}\sum_{n=1}^\infty \lambda_n^{-\gamma }|z|^{-\gamma\alpha}
    \leq c |z|^{(1-\gamma)\alpha-1} \sum_{n=1}^\infty n^{-2\gamma}<c|z|^{(1-\gamma)\alpha-1}.
\end{align*}
Since $E_{\alpha,\alpha}(-\lambda_nz^\alpha)$ is analytic in $z\in \Sigma$, we deduce that the
series in analytic in $t$. These discussions show assertion (i). Next, for any $t>0$, both
series converge uniformly in $[t,\infty)$, and there holds
\begin{align*}
  |e^{-t z}h_1(t)|\leq ce^{-t\Re(z)} t^{-\alpha}
  \Big(\sum_{n=1}^\infty|(u_0,\varphi_n)|^2+\lambda_n^{-2}|(f,\varphi_n)|^2\Big)^\frac12\Big(\sum_{n=1}^\infty n^{-4}\Big)^\frac12\leq ce^{-t\Re(z)}t^{-\alpha},\quad t>0,
\end{align*}
and the function $e^{-t\Re(z)}t^{-\alpha}$ is integrable in $t$ over $(0,\infty)$ for any fixed $z$ with $\Re(z)>0$.
By Lebesgue dominated convergence theorem, we can take Laplace transform termwise and by lemma \ref{lem:mlf}(ii), we obtain
\begin{equation*}
  \int_0^\infty e^{-zt}\sum_{n=1}^\infty \rho_nE_{\alpha,1}(-\lambda_nt^\alpha) \d t = \sum_{n=1}^\infty\rho_n\int_0^\infty e^{-zt}E_{\alpha,1}(-\lambda_nt^\alpha) \d t= \sum_{n=1}^\infty \frac{\rho_nz^{\alpha-1}}{z^\alpha+\lambda_n}, \quad \Re(z)>0.
\end{equation*}
Thus, the Laplace transform of $h_1(t)$ exists.
The Laplace transform of $h_2(t)$ follows similarly from the estimate
\begin{equation*}
|e^{-t z}h_2(t)|\leq ce^{-t\Re(z)} t^{(1-\gamma)\alpha-1} \sum_{n=1}^\infty n^{-2\gamma}
\leq ce^{-t\Re(z)}t^{(1-\gamma)\alpha-1},
\end{equation*} and since $\gamma<1$, the function
$e^{-t\Re(z)}t^{(1-\gamma)\alpha-1}$ is integrable in $t$ over $(0,\infty)$ for any fixed $z$ with $\Re(z)>0$.
Then by Lebesgue dominated convergence theorem and lemma \ref{lem:mlf}(ii), we obtain for $ \Re(z)>0$
\begin{align*}
  \int_0^\infty e^{-zt}\sum_{n=1}^\infty |\varphi_n(0)|^2  t^{\alpha-1}E_{\alpha,\alpha}(-\lambda_nt^\alpha) \d t
   &= \sum_{n=1}^\infty |\varphi_n(0)|^2 \int_0^\infty e^{-zt}
  t^{\alpha-1}E_{\alpha,\alpha}(-\lambda_nt^\alpha) \d t
  =\sum_{n=1}^\infty \frac {\varphi_n(0)^2}{z^\alpha+\lambda_n}.
\end{align*}
This shows assertion (ii), and completes the proof of the proposition.
\end{proof}

\section{Uniqueness}\label{sec:unique}

In this section, we study the uniqueness of the inverse problem: given the observation $h(t)
=u(0,t)$ at the left end point $x=0$, can we uniquely determine the potential $q$ and the order
$\alpha$? Note that without any restriction on the boundary data $g$, the desired uniqueness
result does not hold. This is illustrated by the following example with a zero excitation $g\equiv0$.
\begin{example}\label{exam:counter-unique}
Let $g\equiv0$, $a\equiv1$, and $\alpha\in(0,1)$.
Then consider the following two sets of problem data:
\begin{itemize}
\item[{\rm(a)}] $q\equiv 0$, $f(x)=\frac{\pi^2}{8}(\cos \frac{\pi x}{2}+9\cos\frac{3\pi x}{2})$, $u_0(x)=\frac12\cos\frac{\pi x}{2}+\frac{3}{2}\cos \frac{3\pi x}{2}$;
\item[{\rm(b)}] $\tilde q\equiv \frac{5\pi^2}{4}$, $\tilde f(x)=\frac{9\pi^2}{4}\cos\frac{\pi x}{2} $, $\tilde u_0(x)=2\cos\frac{\pi x}{2}$.
\end{itemize}
Then the eigenvalues $\{\lambda_n(q)\}_{n=1}^\infty$  and $L^2(\Omega)$ orthonormal eigenfunctions $\{\varphi_n\}_{n=1}^\infty$ are given respectively by
\begin{align*}
  \lambda_n(q)=(n-\tfrac{1}{2})^2\pi^2+q\quad \mbox{and}\quad {\varphi_n(x)} = \sqrt{2}\cos((n-\tfrac{1}{2})\pi x),\quad n=1,2,\ldots.
\end{align*}
By Corollary \ref{cor:sol-rep}, the solution $u$ to the direct problem \eqref{eqn:fde} is given by
\begin{equation*}
  u(x,t) = \sum_{n=1}^\infty ((u_0,\varphi_n)-\lambda_n(q)^{-1}(f,\varphi_n))E_{\alpha,1}(-\lambda_n(q)t^\alpha)\varphi_n(x) + \sum_{n=1}^\infty \lambda_n(q)^{-1}(f,\varphi_n)\varphi_n(x).
\end{equation*}
Thus the solutions $u$ and $\tilde u$ for cases {\rm(}a{\rm)} and {\rm(}b{\rm)} are given respectively by
\begin{align*}
  u(x,t)  = \tfrac{1}{2}\cos \tfrac{\pi}{2}x + \tfrac{1}{2}\cos \tfrac{3\pi }{2}x + E_{\alpha,1}(-\tfrac{9\pi^2}{4}t^\alpha)\cos\tfrac{3\pi}{2}x\quad \mbox{and}\quad
  \tilde u(x,t)  =  [1+E_{\alpha,1}(-\tfrac{9\pi^2}{4}t^\alpha)]\cos\tfrac{\pi x}{2}.
\end{align*}
Thus, in both cases, the boundary observation $h$ is given by $h(t)=1+E_{\alpha,1}
(-\tfrac{9\pi^2}{4}t^\alpha),$ and it is impossible to determine the potential
$q$ uniquely from $h$ for $t\in[0,T]$. This shows the impossibility of uniquely
recovering the potential $q$ in the operator $\mathcal{A}$ generally, even if
it {is} assumed to be a constant, and consequently, the desired identifiability result does not hold.
Thus, we excite the system by a nonzero boundary condition $g$ in order to ensure that the data
$h$ contains sufficient information to determine  $q$ uniquely.
\end{example}

Now we proceed to the uniqueness. The proof is split into two steps, and both steps
rely on the time analyticity result in proposition \ref{prop:analyticity}. The first step is
concerned with the unique recovery of the order $\alpha$ and partial information of
the unknown initial data $u_0$ / source $f$. The notation $\mathbb{K}$ denotes the set $\{k\in\mathbb{N}:\rho_k\neq0\}$,
i.e., the support of the sequence $(\rho_0,\rho_1,\ldots)$, with the constants $\rho_k$ defined in
corollary \ref{cor:sol-rep}, and the set $\tilde{\mathbb{K}}$ is defined similarly. It is worth noting that
the set $\mathbb{K}$ is not \textit{a priori} known, since the elliptic operator $\mathcal{A}$ is not
fully known (due to the unknown potential $q$). The condition $\mathbb{K}\neq \emptyset$ holds as long
as $h(t)\not\equiv\mbox{constant}$ on $[0,T_2]$, and thus it is very mild.
\begin{theorem}\label{thm:init}
Let $(q,f,u_0),(\tilde q,\tilde f,\tilde u_0)\in \mathcal{A}\times L^2(\Omega)\times L^2(\Omega)$,
and $h,\tilde h$ be the corresponding observations. Let $0\leq T_1<T_2<\infty$, and the
boundary excitation $g=\tilde g=0$ for $t\in[0,T_2]$. Then the identity $ h(t)=\tilde h(t)$,
$t\in [T_1,T_2]$ implies $\rho_0=\hat \rho_0$, $\{(\rho_k,\lambda_k)\}_{k\in\mathbb{K}}=\{(\tilde \rho_k,\tilde
\lambda_k)\}_{k\in\tilde{\mathbb{K}}}$ and $ \alpha = \tilde\alpha$, if $\mathbb{K}\neq \emptyset$.
\end{theorem}
\begin{proof}
Since $g\equiv0$ for $t\in[0,T_2]$, it follows from Corollary \ref{cor:sol-rep} that $h(t)$
admits a Dirichlet representation
\begin{equation*}
   h(t) =  \rho_0 + \sum_{k\in\mathbb{K}} \rho_kE_{\alpha,1}(-\lambda_kt^\alpha).
\end{equation*}
By proposition \ref{prop:analyticity}(i), $h(t)$ is an analytic function in $t>0$. By analytic continuation, the
condition $h(t)=\tilde h(t)$ for $t\in[T_1,T_2]$ holds implies $h(t)=\tilde h(t)$ for all $t>0$, i.e.
\begin{equation*}
  \rho_0 + \sum_{k\in\mathbb{K}} \rho_kE_{\alpha,1}(-\lambda_kt^\alpha) = \tilde\rho_0+ \sum_{k\in\tilde{\mathbb{K}}} \tilde\rho_k E_{\tilde \alpha,1}(-\tilde \lambda_kt^{\tilde \alpha}).
\end{equation*}
Using the decay property of $E_{\alpha,1}(-\eta)$ in lemma \ref{lem:mlf}(i), cf. \eqref{eqn:bdd-sum}, we derive $\rho_0 = \tilde\rho_0$, $\lambda_0=\tilde \lambda_0$ and hence
\begin{equation*}
   \sum_{k\in\mathbb{K}} \rho_kE_{\alpha,1}(-\lambda_kt^\alpha) =  \sum_{k\in\tilde{\mathbb{K}}} \tilde\rho_k E_{\tilde \alpha,1}(-\tilde \lambda_kt^{\tilde \alpha}),\quad \forall t>0.
\end{equation*}
Now by proposition \ref{prop:analyticity}(ii),  we obtain
\begin{equation*}
  \sum_{k\in\mathbb{K}} \frac{\rho_k z^{\alpha-1}}{z^\alpha + \lambda_k} = \sum_{k\in\tilde{\mathbb{K}}} \frac{\tilde  \rho_k z^{ \tilde \alpha-1}}{z^{\tilde  \alpha} + \tilde \lambda_k}.
\end{equation*}
From this identity we shall deduce $\alpha = \tilde \alpha $ and $\{(\rho_,\lambda_k)\}_{k\in \mathbb{K}} = \{(\tilde\rho_k,\tilde \lambda_k)\}_{k\in \tilde{\mathbb{K}}}$
First, we prove $\alpha=\tilde \alpha$.
Assuming that $\alpha>\tilde \alpha$, dividing both sides by $z^{\tilde \alpha-1}$ and setting $\eta=z^\alpha$, we have
\begin{equation*}
  \sum_{k\in\mathbb{K}} \frac{\rho_k \eta^{1-\frac{\tilde\alpha}{\alpha}}}{\eta + \lambda_k} = \sum_{k\in\tilde{\mathbb{K}}} \frac{\tilde \rho_k}{\eta^{ \frac{\tilde\alpha}{\alpha}} + \tilde \lambda_k}.
\end{equation*}
Choosing arbitrary $k_0\in \mathbb{K}$ and rearranging terms, we derive
\begin{equation*}
 \rho_{k_0} \eta^{1-\frac{\tilde\alpha}{\alpha}} =
  \Big(\sum_{k\in\tilde{\mathbb{K}}} \frac{\tilde \rho_k}{\eta^{ \frac{\tilde\alpha}{\alpha}} + \tilde \lambda_k}
  - \sum_{k\in\mathbb{K}, k\neq k_0} \frac{\rho_k \eta^{1-\frac{\tilde\alpha}{\alpha}} }{\eta + \lambda_k}\Big)(\eta + \lambda_{k_0})
\end{equation*}
Letting $\eta\rightarrow-\lambda_{k_0}$ and noting that $\alpha>\tilde \alpha$, the right hand side of
the identity tends to zero (noting that all $\tilde\lambda_k$ are all real and positive,
and $\arg((-\lambda_{k_0})^{ \frac{\tilde\alpha}{\alpha}}) = \frac{\tilde\alpha\pi}{\alpha} \in (0,\pi)$, and hence  $\rho_{k_0}=0$,
which contradicts the assumption $k_0 \in \mathbb{K}$. Therefore, we deduce
$\alpha\leq \tilde\alpha$. The identical argument yields $\alpha\geq \tilde\alpha$, so we conclude $\alpha=\tilde\alpha$.
The preceding discussion yields%
\begin{equation}\label{eqn:eig-eq}
  \sum_{k\in\mathbb{K}} \frac{\rho_k }{\eta + \lambda_k} = \sum_{k\in\tilde{\mathbb{K}}} \frac{\tilde \rho_k}{\eta + \tilde \lambda_k}.
\end{equation}
In view of the asymptotics of the eigenvalues $\lambda_n$ and $\tilde \lambda_n$, i.e.,
$\lambda_n=O(n^2)$ and $\tilde \lambda_n=O(n^2)$, both sides of the identity converge uniformly in any
compact subset in $\mathbb{C}\setminus(\{-\lambda_k\}_{k\in \mathbb{K}}\cup \{-\tilde\lambda_k\}_{k\in{\tilde{\mathbb{K}}}})$
and are analytic in $\mathbb{C}\setminus(\{-\lambda_k\}_{k\in \mathbb{K}}\cup \{-\tilde\lambda_k\}_{k\in{\tilde{\mathbb{K}}}})$.
Assume that $\lambda_j\not\in \{\tilde\lambda_k\}_{k\in{\tilde{\mathbb{K}}}}$ for some $j\in \mathbb{K}$. Then
we can choose a small circle $C_j$ centered at $-\lambda_j$ and $\{-\tilde\lambda_k\}_{k\in{\tilde{\mathbb{K}}}}$
is not included in the disk centered at $-\lambda_j$ enclosed by $C_j$. Integrating on $C_j$ and applying the
Cauchy theorem, we obtain $\frac{2\pi i\rho_j}{\lambda_j} = 0$, which contradicts the assumption $\rho_j\neq0$.
Hence, $\lambda_j\in \{\tilde\lambda_k\}_{k\in{\tilde{\mathbb{K}}}}$ for every $j\in\mathbb{K}$. Likewise,
$\tilde \lambda_j\in \{\lambda_k\}_{k\in \mathbb{K}}$ for every $j\in\tilde{\mathbb{K}}$. Consequently, we have
proved $\{\lambda_k\}_{k\in\tilde{\mathbb{K}}} = \{\tilde\lambda_k\}_{k\in\tilde{\mathbb{K}}} $, and from \eqref{eqn:eig-eq}, we obtain
\begin{equation*}
  \sum_{k\in\mathbb{K}} \frac{\rho_k-\tilde \rho_k}{\eta+\lambda_k} = 0,\quad \forall \eta\in \mathbb{C}\setminus\{-\lambda_k\}_{k\in\mathbb{K}}.
\end{equation*}
By selecting $j\in\mathbb{K}$, and integrating over $C_j$, we obtain $\frac{2\pi i(\rho_j-\tilde \rho_j)}{\lambda_j}=0$,
which directly implies $\rho_j=\tilde \rho_j.$ This completes the proof of the theorem.
\end{proof}

\begin{remark}
If $f\equiv0$, and $(u_0,\varphi_n)\neq0$, $n\in \mathbb{N}$, then theorem \ref{thm:init}
implies that the sequence $\{(\lambda_n,(u_0,\varphi_n))\}_{n=1}^\infty$ is uniquely determined by
the lateral Cauchy data on $[T_1,T_2]$. However, this does not imply that $u_0$ is uniquely determined yet, since the
potential $q$ and also the eigenfunctions $\varphi_n$ are still unknown. A similar observation can be made when $u_0\equiv0$ and $f\neq0$.
\end{remark}

The next result gives the unique recovery of the potential.
\begin{theorem}\label{thm:main}
Suppose that $(q,f,u_0),(\tilde q,\tilde f,\tilde u_0)\in \mathcal{A}\times L^2(\Omega)\times L^2(\Omega)$. Fix
$0\leq T_1<T_2<T<\infty$, and suppose that the boundary condition $g\in L^\infty(0,T)$ {satisfies} $g=0$ on $[0,T_2]$ and $g\neq0$ on
$[T_2,T]$. Then the identity $h(t)=\tilde h(t)$, $t\in [T_1,T]$ implies $q=\tilde q$.
\end{theorem}
\begin{proof}
In view of the linearity of problem \eqref{eqn:fde}, we can decompose the data $h(t)$ into
\begin{equation*}
  h(t) = u(0,t;0,f,u_0) + u(0,t;g,0,0),\quad t\in (0,T],
\end{equation*}
with the components $u(0,t;0,f,u_0)$ and $u(0,t;0,f,u_0)$ given by
\begin{align*}
  u(0,t;0,f,u_0) &= \sum_{k\in\mathbb{K}} \rho_kE_{\alpha,1}(-\lambda_kt^\alpha),\\
  u(0,t;g,0,0)  &= \sum_{n=1}^\infty \int_0^t (t-s)^{\alpha-1}E_{\alpha,\alpha}(-\lambda_n(t-s)^{\alpha-1})u(s)\d s|\varphi_n(0)|^2,
\end{align*}
respectively, which solves problem \eqref{eqn:fde} with $g\equiv0$ and $f= u_0\equiv0$, respectively, cf. proposition \ref{prop:solrep}.
According to the boundary excitation $g$, the interval $[0,T]$ can be divided into two subintervals:
$(0,T_2]$ and $[T_2,T_1]$. For $t\in(0,T_2)$, it follows directly from theorem \ref{thm:init} that
$\{(\rho_k,\lambda_k)\}_{k\in\mathbb{K}}=\{(\tilde \rho_k,\tilde\lambda_k)\}_{k\in\tilde{\mathbb{K}}}$ and
$\alpha=\tilde\alpha$, from which we have
$u(0,t;0,f,u_0) = \tilde u(0,t;0,\tilde f,\tilde u_0)$ for all $ t>0$.
For $t\in[T_2,T]$, this and the identity $ h(t)=\tilde h(t)$ lead to $u(0,t;g,0,0)=\tilde u(0,t;\tilde g,0,0)$.
This leads to
\begin{align*}
   &\sum_{n=1}^\infty \int_{T_2}^t (t-s)^{\alpha-1}E_{\alpha,\alpha}(-\lambda_n(t-s)^{\alpha-1})g(s)\d s|\varphi_n(0)|^2\\
  =&\sum_{n=1}^\infty \int_{T_2}^t (t-s)^{\alpha-1}E_{\alpha,\alpha}(-\tilde \lambda_n(t-s)^{\alpha-1})g(s)\d s|\tilde\varphi_n(0)|^2,\quad t\in[T_2,T]
\end{align*}
Since $g\in L^2(T_2,T)$ is nonzero for almost all $t\in (T_2,T)$ and the kernel belongs to $L^1(0,\infty)$, it
follows from Titchmarsh convolution theorem \cite[Theorem VII]{Titchmarsh:1926} that
\begin{equation*}
  \sum_{n=1}^\infty t^{\alpha-1}E_{\alpha,\alpha}(-\lambda_nt^{\alpha-1})
  |\varphi_n(0)|^2= \sum_{n=1}^\infty t^{\alpha-1}E_{\alpha,\alpha}(-\tilde\lambda_nt^{\alpha-1})|\tilde\varphi_n(0)|^2,\quad t\in [0,T-T_2].
\end{equation*}
Now by the analyticity of the functions on both sides, cf. proposition \ref{prop:analyticity}(i), we have
\begin{equation*}
  \sum_{n=1}^\infty t^{\alpha-1}E_{\alpha,\alpha}(-\lambda_nt^{\alpha-1})
  |\varphi_n(0)|^2= \sum_{n=1}^\infty t^{\alpha-1}E_{\alpha,\alpha}(-\tilde\lambda_nt^{\alpha-1})|\tilde\varphi_n(0)|^2,\quad t\in (0,\infty).
\end{equation*}
Note that $\varphi_n(0)\neq0$ for all $n\in\mathbb{N}$, similar to the proof of theorem \ref{thm:init}, one
can show
\begin{equation*}
  \{(\lambda_n,|\varphi_n(0)|)\}_{n\in\mathbb{N}} = \{(\tilde\lambda_n,|\tilde\varphi_n(0)|)\}_{n\in\mathbb{N}}.
\end{equation*}
Finally by the classical Gel'fand-Levitan theory \cite{GelfandLevitan:1955,LevitanSargsjan:1975}, we deduce $q=\tilde q$.
\end{proof}

\begin{corollary}
If one of the functions $u_0$ and $f$ is zero, then the other can be uniquely determined from the observation $h(t)$, $t\in [0,T]$.
\end{corollary}
\begin{proof}
We consider the case $f\equiv0$, and the other case $u_0\equiv0$ follows similarly.
By theorem \ref{thm:main}, the potential $q$ is uniquely determined by $h(t)$, $t\in [T_1,T]$, and thus also
the eigenfunctions $\varphi_n(x)$ associated with the corresponding elliptic operator $\mathcal{A}$. Then by
theorem \ref{thm:init}, the sequence $\{(u_0,\varphi_k)\}_{k\in\mathbb{K}}$
is uniquely determined, which directly gives the unique recovery of the initial data $u_0$.
\end{proof}

\begin{remark}\label{rmk:JingYamamoto}
There have been several works on identifying multiple parameters from one single observation
\cite{XianYanWei:2020,KianLiLiu:2020,JingYamamoto:2021}. The recent work \cite{JingYamamoto:2021}
is closest to the current one in some sense, which is concerned with the following model on
$\Omega=(0,1)$, with $\alpha\in(0,2)$,
\begin{equation*}
  \left\{\begin{aligned}
     \partial_t^\alpha u - \mathcal{A}u  & = 0, \quad\mbox{in }\Omega\times(0,T],\\
     a(0)\partial_xu(0,t)-hu(0,t) & = 0, \quad \mbox{in }(0,T],\\
     a(1)\partial_xu(1,t)+H(1,t) & = 0,\quad \mbox{in } (0,T],\\
     u(0) &= u_0, \quad \mbox{in }\Omega,\\
     u'(0)& = u_0',\quad \mbox{in }\Omega, \mbox{ if } \alpha\in (1,2).
  \end{aligned}\right.
\end{equation*}
The inverse problem is to recover $u_0$, $q$, $\alpha$, $h$ and $H$
from two boundary observations, i.e., $u(0,t)$ and $u(1,t)$. They
proved the uniqueness of the recovery under the following condition
{\rm(}for $\alpha\in (0,1)${\rm)}: $(u_0,\varphi_n)\neq0$, for all $n\in\mathbb{N}$ \cite[Theorem 1]{JingYamamoto:2021}.
This condition assumes that all the eigenmodes of the initial value $u_0$
should be nonzero, which is generally restrictive, and can be relaxed
using multiple initial conditions \cite[Theorem 1']{JingYamamoto:2021}.
In contrast, theorem \ref{thm:main} relies on the nonzero boundary
excitation $g(t)$ for the potential recovery, and thus avoids this assumption.
\end{remark}

\begin{remark}
There are several potential extensions of the stated uniqueness results.
(1) The results hold also for the multi-term time-fractional model, which involves multiple time
fractional derivatives, i.e., the term $\partial_t^\alpha u$ in the model \eqref{eqn:fde} is replaced
by $\sum_{i=1}^Nr_i\partial_t^{\alpha_i} u$, with $r_i>0$ and $0<\alpha_1<\ldots<\alpha_N<1$. Then
the weights $r_i$ and $\alpha_i$ are uniquely determined, provided that $h(0)\neq0$. (2) One can
uniquely determine the diffusion coefficient $a$ when the potential $q$ is known, by a different version of
Gel'fand-Levitan theory \cite{ChengNakagawa:2009}. (3) The boundary conditions can be of more general
Sturm-Liouville form. Then the Robin coefficients in the boundary conditions can also be determined uniquely
from lateral Cauchy data \cite{JingYamamoto:2021}, cf. remark \ref{rmk:JingYamamoto}.
\end{remark}

The preceding analysis indicates that the both steps rely essentially on unique continuation,
which is well known to be severe ill-conditioned. A natural question is how the fractional
paradigm actually affects the degree of ill-conditioning, measured in terms of the asymptotic decay
rate of the singular value spectrum of the associated (linearized) forward map. This issue has been
numerically studied for several inverse problems in \cite{JinRundell:2015}; see also
\cite{RundellYamamoto:2018} for the inverse potential problem. However, a
theoretical analysis in the context of potential recovery from lateral Cauchy data is still unavailable.

\section{Reconstruction algorithm}\label{sec:alg}

Now we describe an algorithm for simultaneously recovering the potential $q$, the order
$\alpha$, and also $u_0$, under the assumption $f\equiv0$ (or also $f$, if $u_0\equiv0$).
The procedure is directly inspired by the uniqueness proof, and consists of two steps.

\subsection{Step 1: order determination and numerical continuation}
In the first step, we determine the fractional order $\alpha$ and numerically continuate
the trace data $h(t)$ from $[0,T_1]$ to the whole interval $[0,T]$ (to assist the recovery
of the potential $q$). We discuss the two issues separately. The recovery of the order $\alpha$
cannot be carried out in the usual manner by means of least-squares fitting, since the problem data in the direct problem
\eqref{eqn:fde} over the interval $[0,T_1]$ is not fully known. The next result suggests
one possible recovery formula for the order $\alpha$ from the small time asymptotics of the
observation $h(t)$, under suitable smoothness condition $u_0$ and $f$, if the function
$\mathcal{A}u_0+f$ does not vanish at $x=0$.
\begin{proposition}\label{prop:asymptotic}
If $u_0\in D(A^{1+s})$ and $f\in D(A^s)$ with $s\in (\frac14,\frac54]$, then for any $\epsilon\in(0,s-\frac14)$,
$h(t)=u(0,t)$ satisfies the following asymptotic
\begin{equation*}
  h(t) = u_0(0)-(\mathcal{A}u_0(0)+f(0))t^\alpha + O({t^{(1-\frac14-\frac\epsilon2+s)\alpha}}),\quad \mbox{as } t\to 0^+.
\end{equation*}
\end{proposition}
\begin{proof}
By the definition of the Mittag-Leffler function $E_{\alpha,1}(z)$, we have
\begin{equation*}
  E_{\alpha,1}(-\lambda_nt^\alpha)=1-\lambda_nt^\alpha + \lambda_n^2t^{2\alpha}E_{\alpha,1+2\alpha}(-\lambda_nt^{\alpha}).
\end{equation*}
This and the solution representation from corollary \ref{cor:sol-rep} lead to
\begin{align*}
  u(x,t) &= \sum_{n=1}^\infty \big([(u_0,\varphi_n) - \lambda_n^{-1}(f,\varphi_n)]E_{\alpha,1}(-\lambda_nt^\alpha)+\lambda_n^{-1}(f,\varphi_n)\big)\varphi_n(x)\\
  & = \sum_{n=1}^\infty (u_0,\varphi_n) \varphi_n(x) -\sum_{n=1}^\infty \big([(u_0,\varphi_n) - \lambda_n^{-1}(f,\varphi_n)]\lambda_nt^\alpha \varphi_n(x)\\
  &\quad + \sum_{n=1}^\infty [(u_0,\varphi_n)-\lambda_n^{-1}(f,\varphi_n)]\lambda_n^2t^{2\alpha}E_{\alpha,1+2\alpha}(-\lambda_nt^\alpha)\varphi_n(x).
\end{align*}
We denote the last sum by $\rm I$. Since the eigenfunctions $\{\varphi_n\}_{n=1}^\infty$ forms an orthonormal basis in $L^2(\Omega)$, by integration by parts, we have
\begin{equation*}
  \sum_{n=1}^\infty\lambda_n(u_0,\varphi_n)\varphi_n(x) = \sum_{n=1}^\infty (u_0,-\mathcal{A}\varphi_n)\varphi_n(x) = \sum_{n=1}^\infty (-\mathcal{A}u_0,\varphi_n)\varphi_n(x) = -\mathcal{A}u_0(x)
\end{equation*}
and $\sum_{n=1}^\infty (f,\varphi_n)\varphi_n(x) = f(x)$. By lemma \ref{lem:mlf}(i), we bound the sum {\rm I} by
\begin{align*}
 \|{\rm I}\|_{L^2(\Omega)}\leq & ct^{4\alpha}\sum_{n=1}^\infty \lambda_n^{\frac12+\epsilon + 2-2s}[\lambda_n^{2+2s}(u_0,\varphi_n)^2+\lambda_n^{2s}(f,\varphi_n)^2]E_{\alpha,1+2\alpha}(-\lambda_nt^{\alpha})^2\\
  \leq & ct^{(2-\frac12-\epsilon+2s)\alpha}\sum_{n=1}^\infty \frac{(\lambda_nt^\alpha)^{\frac12+\epsilon + 2-2s}}{(1+\lambda_nt^\alpha)^2}[\lambda_n^{2+2s}(u_0,\varphi_n)^2+\lambda_n^{2s}(f,\varphi_n)^2]\\
  \leq & ct^{(2-\frac12-\epsilon+2s)\alpha}\sum_{n=1}^\infty[\lambda_n^{2+2s}(u_0,\varphi_n)^2+\lambda_n^{2s}(f,\varphi_n)^2]\\
  =&ct^{(2-\frac12-\epsilon+2s)\alpha}[\|A^{1+s}u_0\|_{L^2(\Omega)}^2+\|A^sf\|_{L^2(\Omega)}],
\end{align*}
where the last inequality follows from the conditions $s\in (\frac14,\frac54]$ and $\epsilon\in(0,s-\frac14)$.
Combining the preceding estimates with the Sobolev embedding theorem directly shows the assertion.
\end{proof}

By proposition \ref{prop:asymptotic}, under mild conditions, the trace data $h(t)$ satisfies
\begin{equation*}
  h(t) = c_0 + c_1 t^\alpha + o(t^{\alpha}),\quad \mbox{as }t\to 0^+.
\end{equation*}
This motivates a simple procedure: minimize over $\alpha$ (and $c_0$ and $c_1$)
the following objective
\begin{equation*}
  J(\alpha,c_0,c_1) = \tfrac12\|c_0+c_1t^\alpha-h(t)\|_{L^2(0,t_0)}^2,
\end{equation*}
for some $t_0>0$ sufficiently close to zero. The minimization can be carried out
by any stand-alone algorithms, e.g., gradient descent, and Newton method.
Note that it is important to take $t_0$ sufficiently close to zero so that
the term $o(t^\alpha)$ is indeed negligible.

Next, we numerically continuate the given data $h(t)$ from the interval
$[0,T_1]$ to $[T_1,T]$, in order to extract the combined information on $u_0$
and $f$. Mathematically, this amounts to recovering $\{(\rho_k,\lambda_k)\}_{k\in\mathbb{K}}$
(with the index set $\mathbb{K}$ defined in section \ref{sec:unique}) from $h$:
\begin{equation*}
    h(t) = \sum_{k\in\mathbb{K}}\rho_kE_{\alpha,1}(-\lambda_kt^\alpha),\quad t\in [0,T_1].
\end{equation*}
By theorem \ref{thm:init}, $\{(\rho_k,\lambda_k)\}_{k\in\mathbb{K}}$ can indeed be uniquely determined
by $h(t)$, $t\in[0,T_1]$. This problem is also known as an infinite-dimensional spectral
estimation problem for $\alpha=1$, for which the issue is to recover $(\rho_k,\lambda_k)$ of
an exponential family \cite{AvdoninGesztesy:2010} and there are several efficient methods for recovery,
e.g., matrix pencil method \cite{HuaSarkar:1990} and MUSIC (MUltiple SIgnal Classification) \cite{Schmidt:1986}.
However, for $\alpha\neq 1$, to the best of our knowledge, there is no known analogue of these methods.
This is essentially due to the inequality $E_{\alpha,1}(-t_1^\alpha)E_{\alpha,2}(-t_2^\alpha) \neq E_{\alpha,1}
(-(t_1+t_2)^\alpha)$ for $\alpha\in(0,1)$, $t_1,t_2\in(0,\infty)$. Instead, we resort to the classical rational
approximation for numerical continuation, i.e.,
\begin{equation*}
  h(t) \approx \frac{p_0+p_1t+\ldots +p_rt^r}{q_0+q_1t + \ldots q_rt^r}:= h_r(t),\quad t\in[0,T_1],
\end{equation*}
where $r\in\mathbb{N}$ is the polynomial order.
The approximation $h_r(t)$ can be constructed efficiently when $h(t)$ is accurate using the
AAA algorithm \cite{NakatsukasaSeteTrefethen:2018}, despite the well-known ill-posed nature of analytic
continuation. This choice is in part motivated by the fact that the function $E_{\alpha,1}
(-\lambda t^\alpha)$ admits excellent rational approximations \cite[Theorem 3.6]{Jin:2021}.
Our numerical experiments indicate that the procedure is indeed viable for exact data.

\subsection{Step 2: recovering $q$ (and $u_0$) by iterative regularization}
With the analytic continuation in Step 1, we can proceed to the reconstruction of the
potential $q$, as in the proof of theorem \ref{thm:main}. Specifically, let
$$
\bar h(t) = \left\{\begin{aligned}
   0, &\quad t\in [0,T_1],\\
   h(t)- h_r(t), &\quad t\in[T_1,T],
\end{aligned}\right.
$$
which represents the reduced data for the boundary excitation $g$ only
(supported on the interval $[T_1,T]$, by construction). This naturally motivates
approximately minimizing
\begin{equation}\label{eqn:Tikh}
  J(q) :=\tfrac12\|F(q)-\bar h\|_{L^2(T_1,T)}^2,
\end{equation}
with $F(q)=u(q)(0,t)$, where $u(q)$ denotes the solution to the direct problem \eqref{eqn:fde}
corresponding to the elliptic operator $\mathcal{A}$, with $u_0\equiv f\equiv 0$
and given $g$. The map $F$ is nonlinear, and one may apply
standard iterative regularization methods \cite{EnglHankeNeubauer:1996,ItoJin:2015},
e.g., (nonlinear) conjugate gradient method. In the numerical experiments, we
employ the conjugate gradient method \cite{AlifanovArtyukhin:1995}, which generally
enjoys fast convergence.

Once the potential $q$ is determined, one can also attempt recovering the initial data $u_0$
from the observation $h(t)$ over the interval $[0,T_1]$, if $f\equiv0$. This can be achieved
by approximately minimizing
\begin{equation*}
  J(u_0) = \tfrac{1}{2}\|F(u_0)-h(t)\|_{L^2(0,T_1)}^2,
\end{equation*}
with $F(u_0)=u(0,t)$, where $u$ denotes the solution to the direct problem \eqref{eqn:fde}
corresponding to the elliptic operator $\mathcal{A}$ with the recovered $q$, and $ f\equiv 0$
(over the interval $(0,T_1)$). The optimization can be carried out efficiently by standard
gradient type methods.

In practice, the gradients of the functionals $J(q)$ and $J(u_0)$ can be
computed efficiently by the adjoint technique. We provide relevant details in the appendix.

\section{Numerical results and discussions}\label{sec:numer}

Now we present some numerical results to illustrate the feasibility of simultaneously recovering the
coefficient $q$ and the fractional order $\alpha$, without fully knowing the
direct problem \eqref{eqn:fde}. The domain $\Omega$ is taken to be the unit interval $\Omega=
[0,1]$, and the final time $T=1$, and $T_1=0.5$. The direct and adjoint problems
are all discretized by the standard continuous piecewise linear Galerkin method in space, and
backward Euler convolution quadrature in time \cite{JinLazarovZhou:2016sisc}. 
The domain $\Omega$ is divided into $M$ subintervals each of width $1/M$. For
the inversion step, we take $M=200$ and $N=2000$. The exact data $h^\dag$ on the lateral boundary
$(0,T)$ is obtained by solving the direct problem \eqref{eqn:fde} on a finer mesh. It is known that
due to the severe ill-conditioning of the inverse problem, the numerical recovery in the presence
of data noise is very challenging. Indeed, for the inverse potential problem, it was observed
numerically in \cite{RundellYamamoto:2018} that there are only a few significant singular values,
and this is also partly confirmed by \cite{SunWei:2017}. This is further
complicated by the unknown problem data in the present context. Thus, our experiments below focus on exact
data. We illustrate on the following two settings, with $u_0$ and $f$ being unknown:
\begin{itemize}
  \item[(i)] $a\equiv1$, $q^\dag = x(1-x)$, $u_0 = x^2(1-x)+\cos(\frac\pi2 x)$ and $f\equiv0$, and $g=\chi_{[T_1,T]}$;
  \item[(ii)] $a\equiv1$, $q^\dag = \min(x,1-x)$, $u_0 = \cos(\frac{3\pi}{2} x)$ and $f\equiv0$, and $g=\chi_{[T_1,T]}$,
\end{itemize}
where $\chi_S$ denotes the characteristic function of the set $S$. The initial condition $u_0$ is
taken to be in $D(A)$ so that the asymptotic expansion in proposition \ref{prop:asymptotic} is
indeed valid. Case (i) involves a smooth potential, and case (ii) a nonsmooth potential.

First we study the recovery of the fractional order $\alpha$ using a least-squares fitting as
described in section \ref{sec:alg}. This procedure relies on the validity of the asymptotic
expansion in proposition \ref{prop:asymptotic}. The recovered orders are presented in
Table \ref{tab:order}, where the minimization is carried out by the L-BFGS-B \cite{ByrdLuNocedal:1995}, with
the box constraint $\alpha \in[0,1]$, using the public implementation \url{https://ww2.mathworks.cn/matlabcentral/fileexchange/35104-lbfgsb-l-bfgs-b-mex-wrapper}
(last accessed on May 20, 2021). Note that the least-squares functional is fraught with many
local minimum, and a good initial guess for $\alpha$ is needed in order to recover the correct order.
It is observed that the accuracy of the recovery tends to
improve as the interval $(0,t_0)$ used in the least-squares formulation shrinks, since
the model function in proposition \ref{prop:asymptotic} represents an increasingly better
approximation as $t\to0^+$. Further, as $\alpha$ increases, the
size of the interval $(0,t_0)$ can be increased without sacrificing the accuracy of the
recovery since the asymptotic expansion is then valid in a larger neighborhood. Thus one
may conclude that with the interval $(0,t_0)$ chosen properly (and of course only for very
accurate data), the order $\alpha$ can indeed be recovered reliably by the least-squares fitting. These
observations hold for both cases (i) and (ii), and thus the smoothness of the potential $q$ does not
seem to influence much the recovery of the order $\alpha$.

\begin{table}[hbt!]
  \centering
  \caption{The recovered order $\alpha$ based on least-squares fitting.\label{tab:order}}
  \begin{threeparttable}
  \subfloat[case (i)]{
  \begin{tabular}{ccccc}
  \toprule
  $t_0\backslash\alpha$ & $0.3$ & $0.5000$ & $0.7000$ & $0.9000$\\
  \midrule
  1e-3 &  0.2488  &  0.6239  &  0.8641  &  1.0000\\
  1e-4 &  0.3208  &  0.6256  &  0.8110  &  0.9647\\
  1e-5 &  0.3631  &  0.6027  &  0.7626  &  0.9000\\
  1e-6 &  0.3760  &  0.5747  &  0.7315  &  0.9000\\
  1e-7 &  0.3737  &  0.5495  &  0.7000  &  0.9000\\
  1e-8 &  0.3665  &  0.5306  &  0.7000  &  0.9000\\
  1e-9 &  0.3570  &  0.5000  &  0.7000  &  0.9000\\
  1e-10&  0.3468  &  0.5000  &  0.7000  &  0.9000\\
  \bottomrule
  \end{tabular}}\quad
  \subfloat[case (ii)]{
  \begin{tabular}{ccccc}
  \toprule
  $t_0\backslash\alpha$ & $0.3$ & $0.5000$ & $0.7000$ & $0.9000$\\
  \midrule
  1e-3 &  0.0008  &  0.2723  &  0.6274  &  0.8829\\
  1e-4 &  0.0271  &  0.4126  &  0.6850  &  0.8980\\
  1e-5 &  0.1214  &  0.4700  &  0.6969  &  0.9000\\
  1e-6 &  0.1932  &  0.4897  &  0.6992  &  0.9000\\
  1e-7 &  0.2398  &  0.4960  &  0.7000  &  0.9000\\
  1e-8 &  0.2667  &  0.4980  &  0.7000  &  0.9000\\
  1e-9 &  0.2813  &  0.5000  &  0.7000  &  0.9000\\
  1e-10&  0.2888  &  0.5000  &  0.7000  &  0.9000\\
  \bottomrule
  \end{tabular}}
\end{threeparttable}
\end{table}

\begin{figure}[hbt!]
\begin{tabular}{cc}
 \includegraphics[width=.48\textwidth]{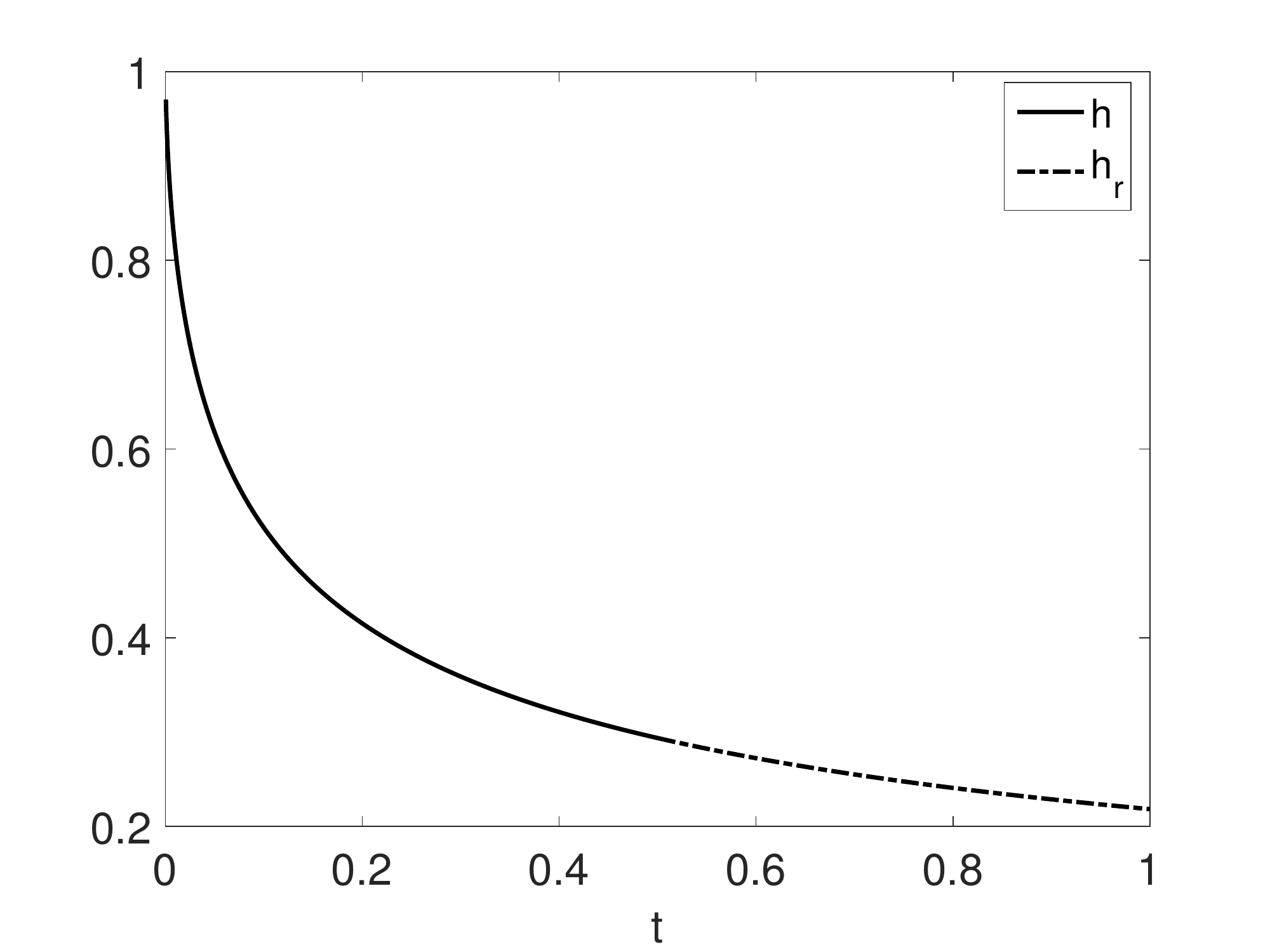} & \includegraphics[width=.48\textwidth]{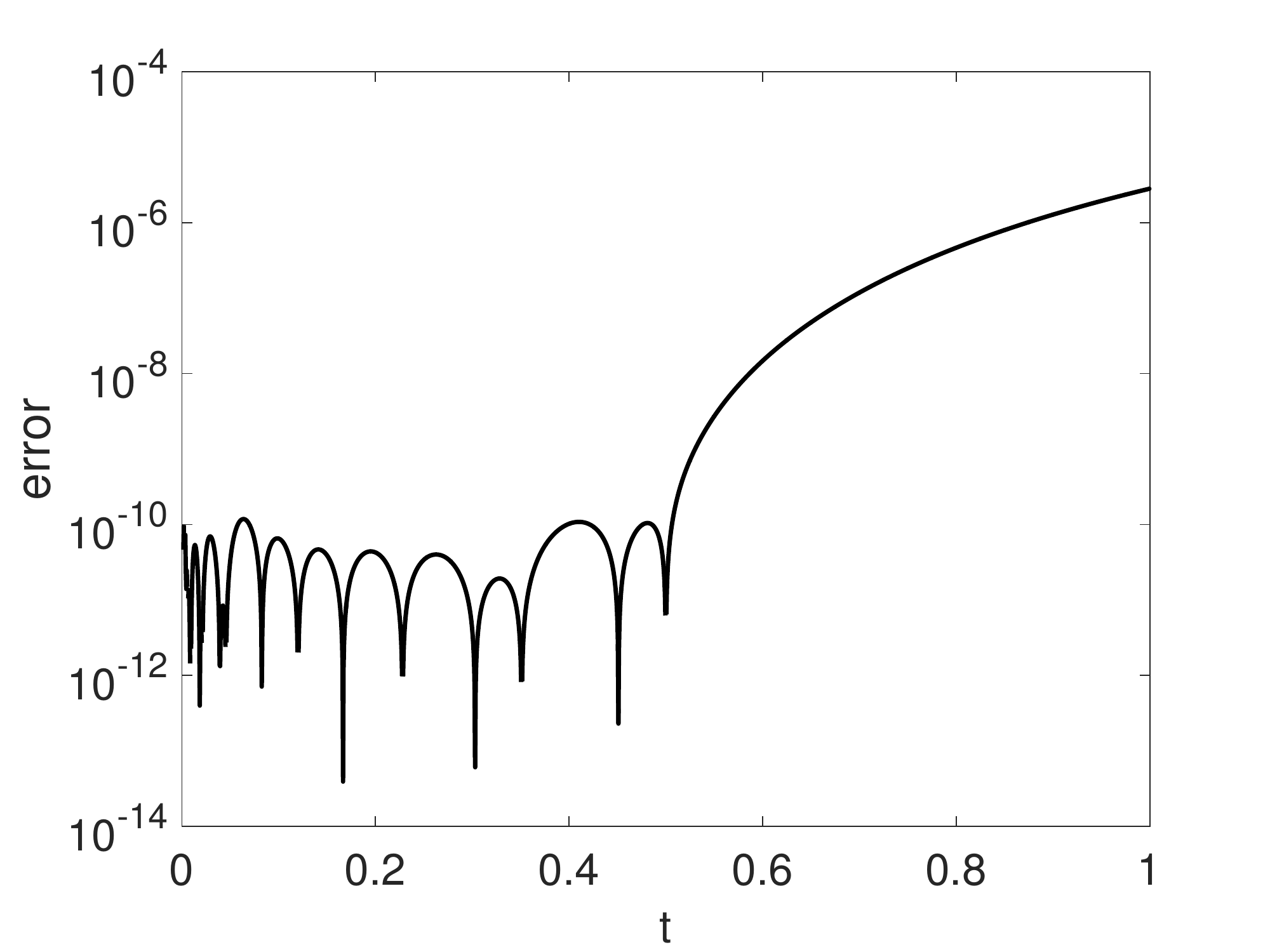}\\
\includegraphics[width=.48\textwidth]{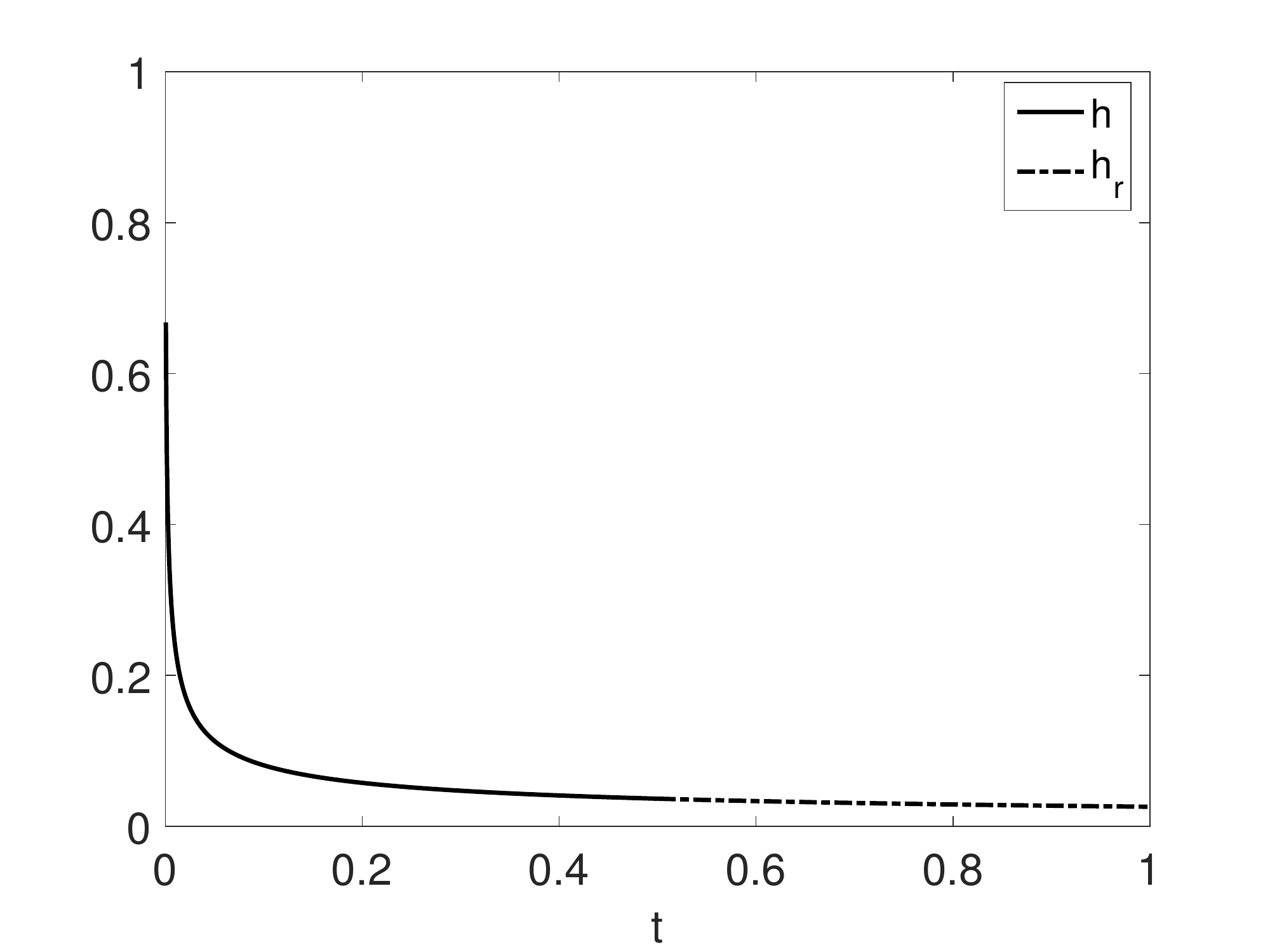} & \includegraphics[width=.48\textwidth]{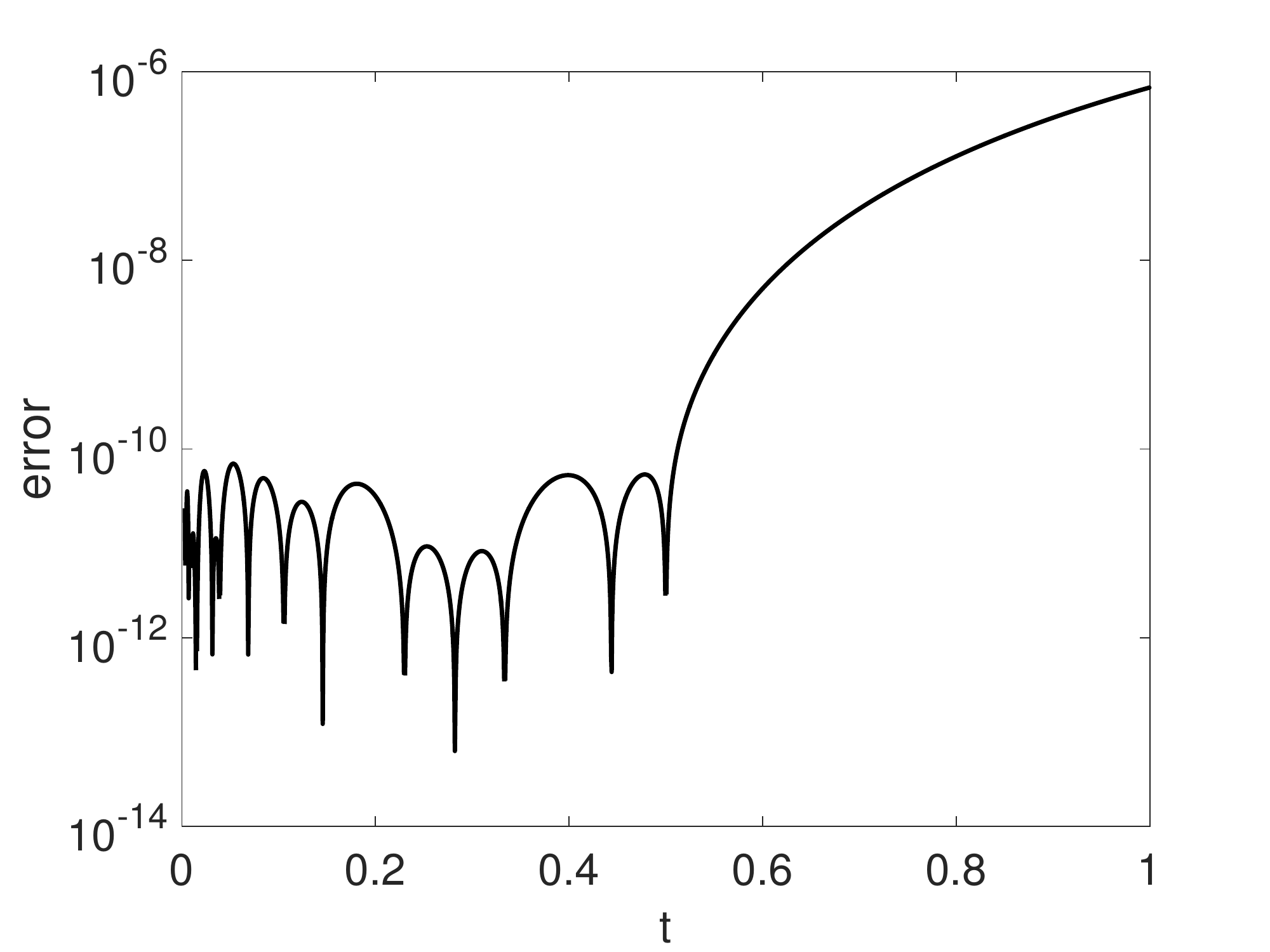}\\
(a) $h^*$ and $h_r$ & (b) error
\end{tabular}
\caption{The analytic continuation $h_r$ of the true data $h^*$ from $[0,T_1]$ to $[T_1,T]$
by rational approximation, and its pointwise error for cases (i) (top) and (ii) (bottom).\label{fig:ac}}
\end{figure}

One step of the recovery procedure is analytic continuation, extending
the observation data $h$ by a rational model $h_r$ from the interval $[0,T_1]$ to
$[T_1,T]$. This step extracts relevant information from unknown initial condition $u_0$
(and source $f$), and plays a central role in formulating the optimization problem
for recovering the potential $q$. This is illustrated in Fig. \ref{fig:ac} for the two
cases at $\alpha=0.5$, where the rational approximation $h_r$ is constructed by the
AAA algorithm \cite{NakatsukasaSeteTrefethen:2018} using the \texttt{MATLAB}
implementation given therein with a tolerance 1e-9, with the resulting $h_r$ of
degree $r=11$. The pointwise error is evaluated against the ground-truth $h^*$
(i.e., $h^*=u(0,t)$, $t\in[0,T]$, $u$ being the solution of the direct problem
\eqref{eqn:fde} with $g\equiv0$) over the interval $[T_1,T]$. (Numerically, larger tolerances, e.g.,
1e-6, can still give an accurate approximation.) Clearly,
$h_r$ does give a fairly accurate approximation to $h^*$, and
the accuracy degrades as one moves away from the interpolating interval $[0,T_1]$.
It is noted that the continuation results for other cases exhibit very similar
behavior. Thus, the rational approximation is a very effective approach for analytic
continuation when exact data is available.

\begin{figure}[hbt!]\setlength{\tabcolsep}{0pt}
\begin{tabular}{ccc}
 \includegraphics[width=.33\textwidth]{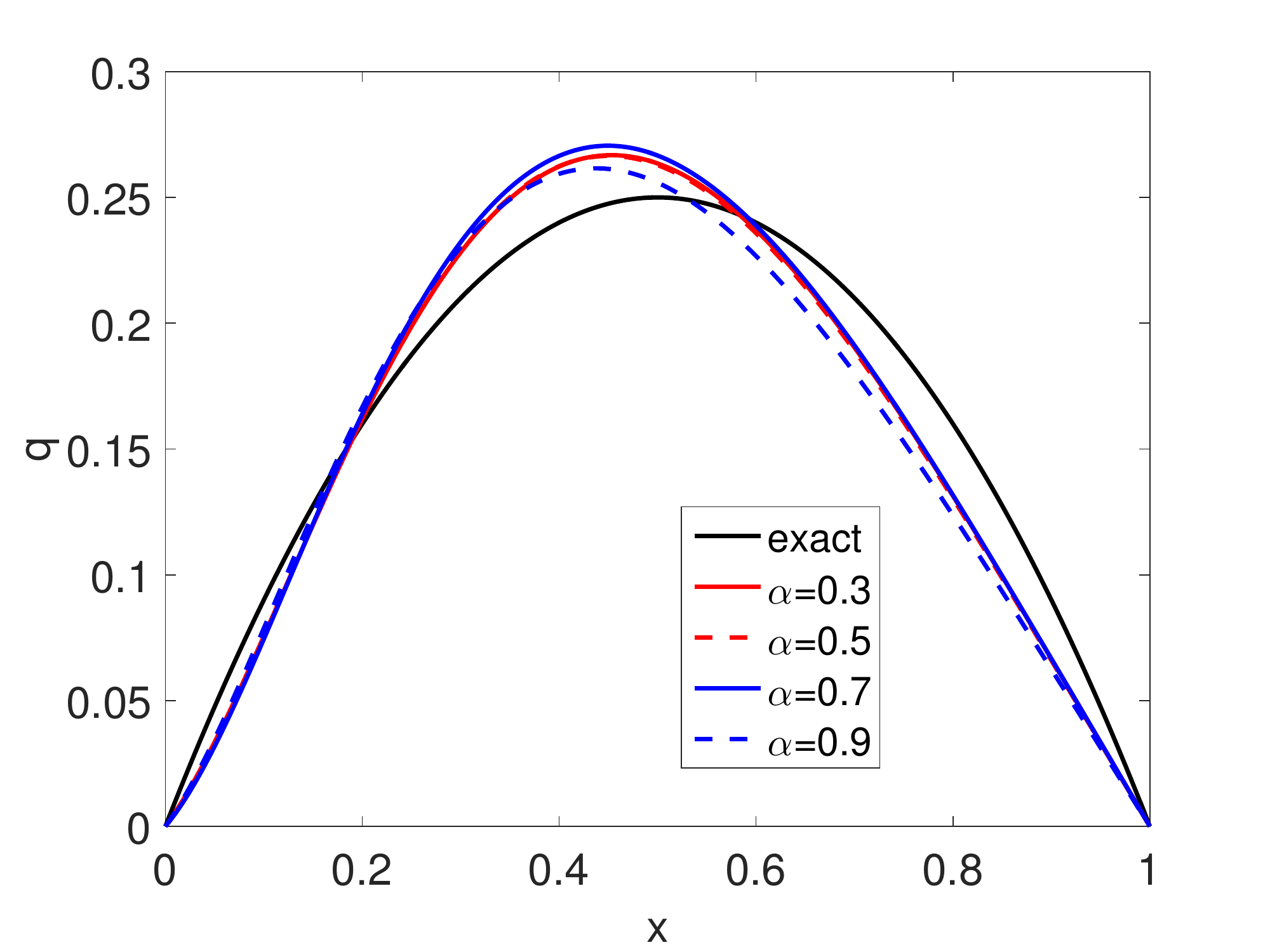} & \includegraphics[width=.33\textwidth]{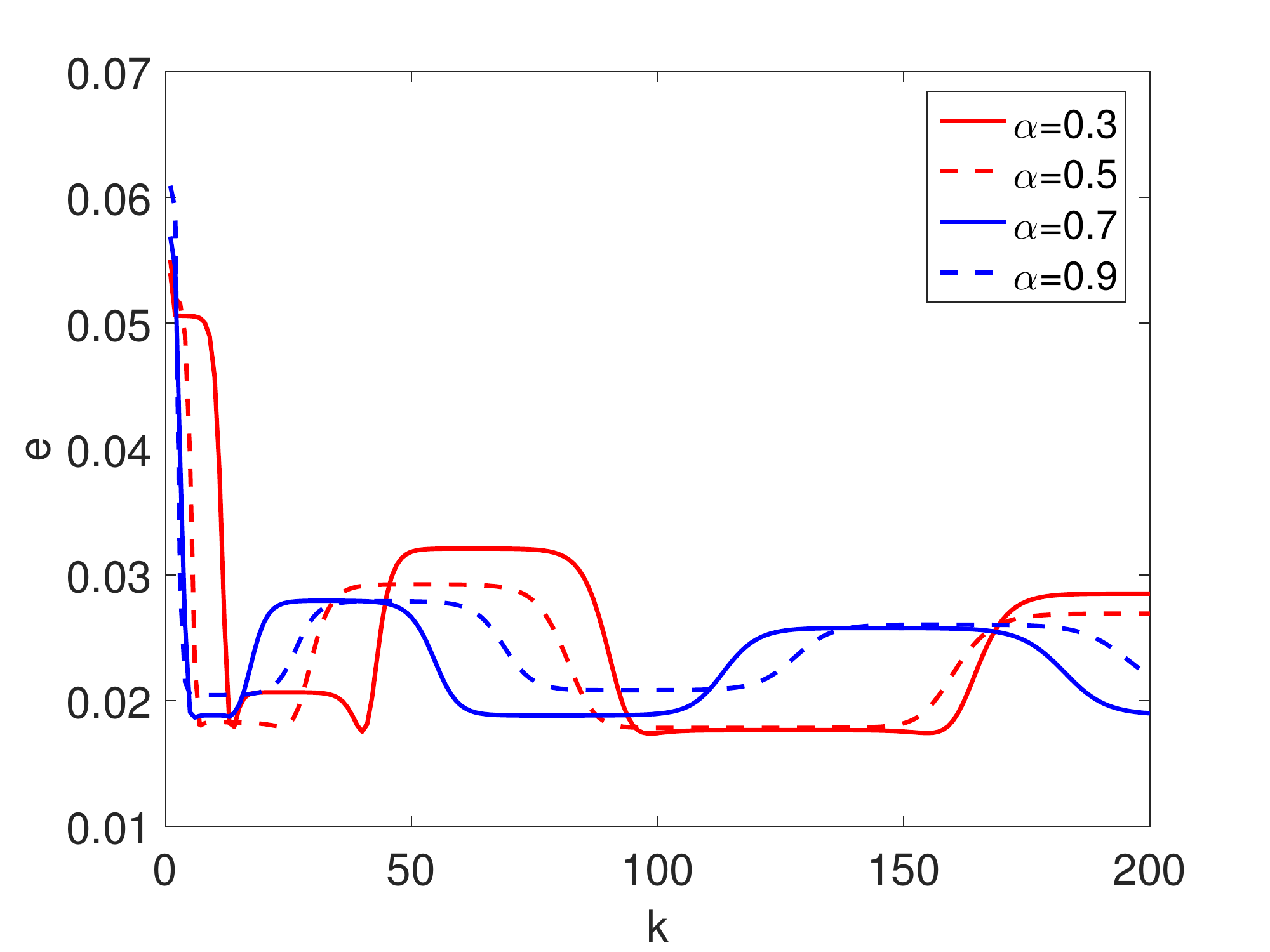} & \includegraphics[width=.33\textwidth]{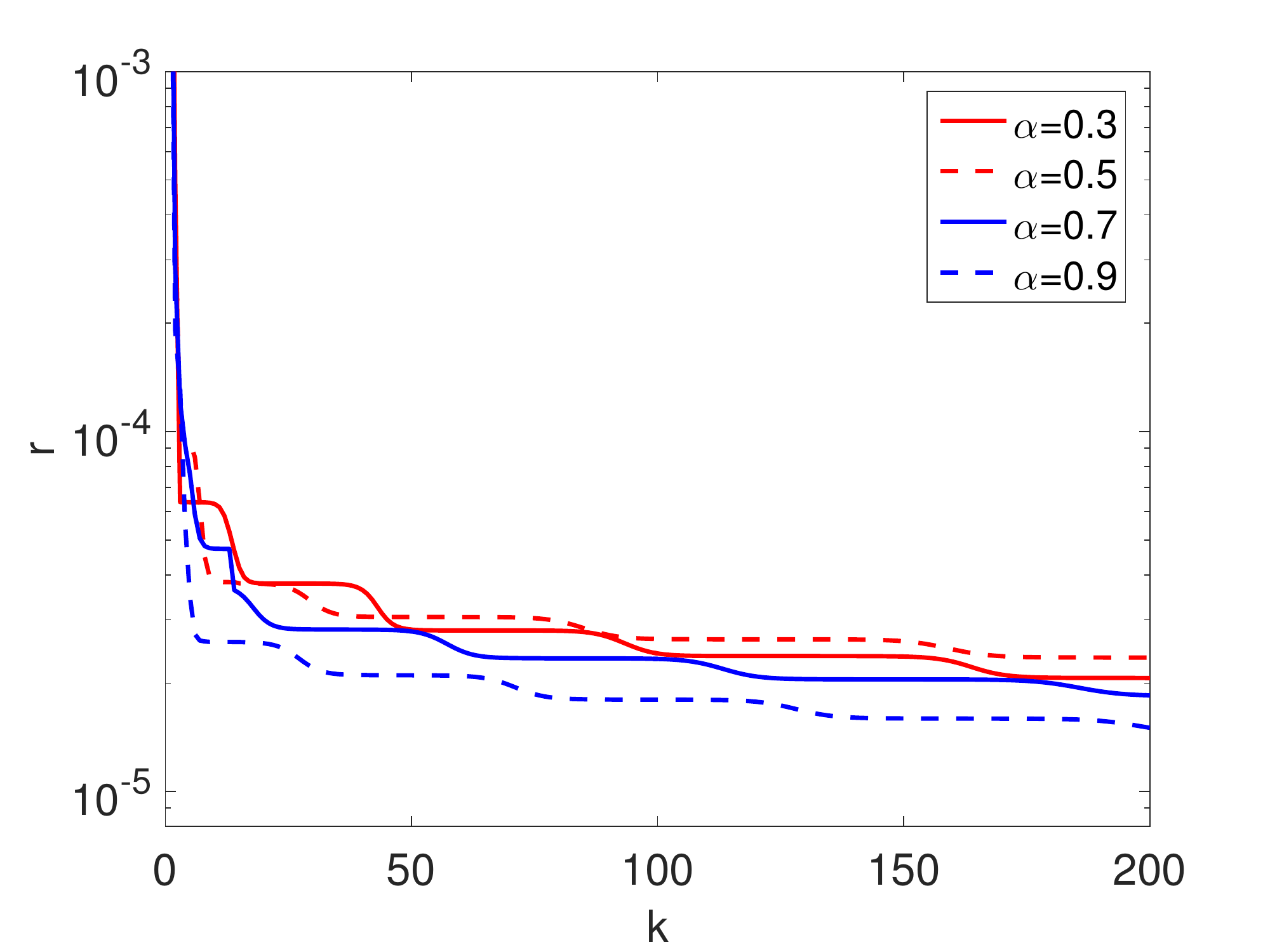}\\
  \includegraphics[width=.33\textwidth]{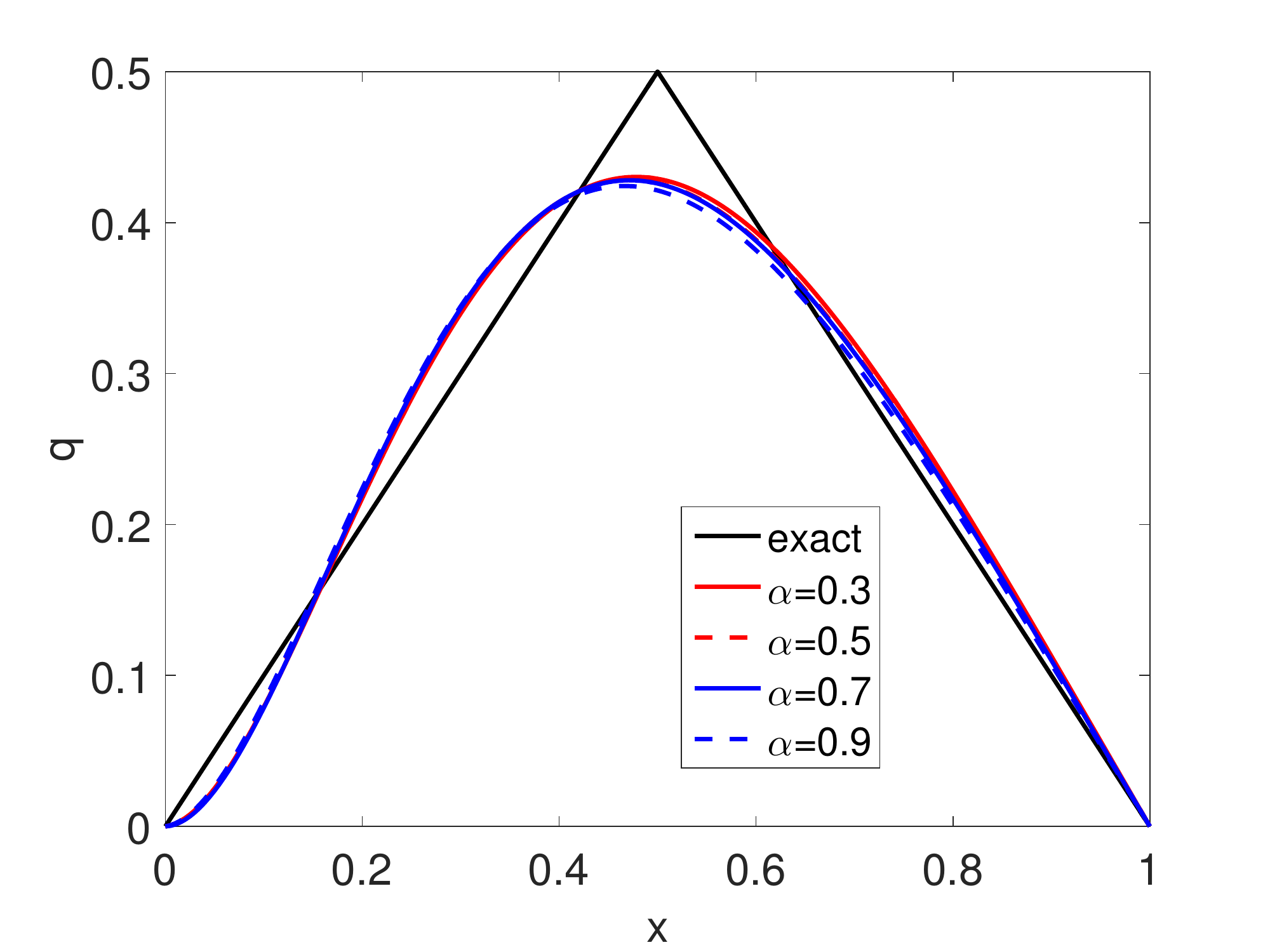} & \includegraphics[width=.33\textwidth]{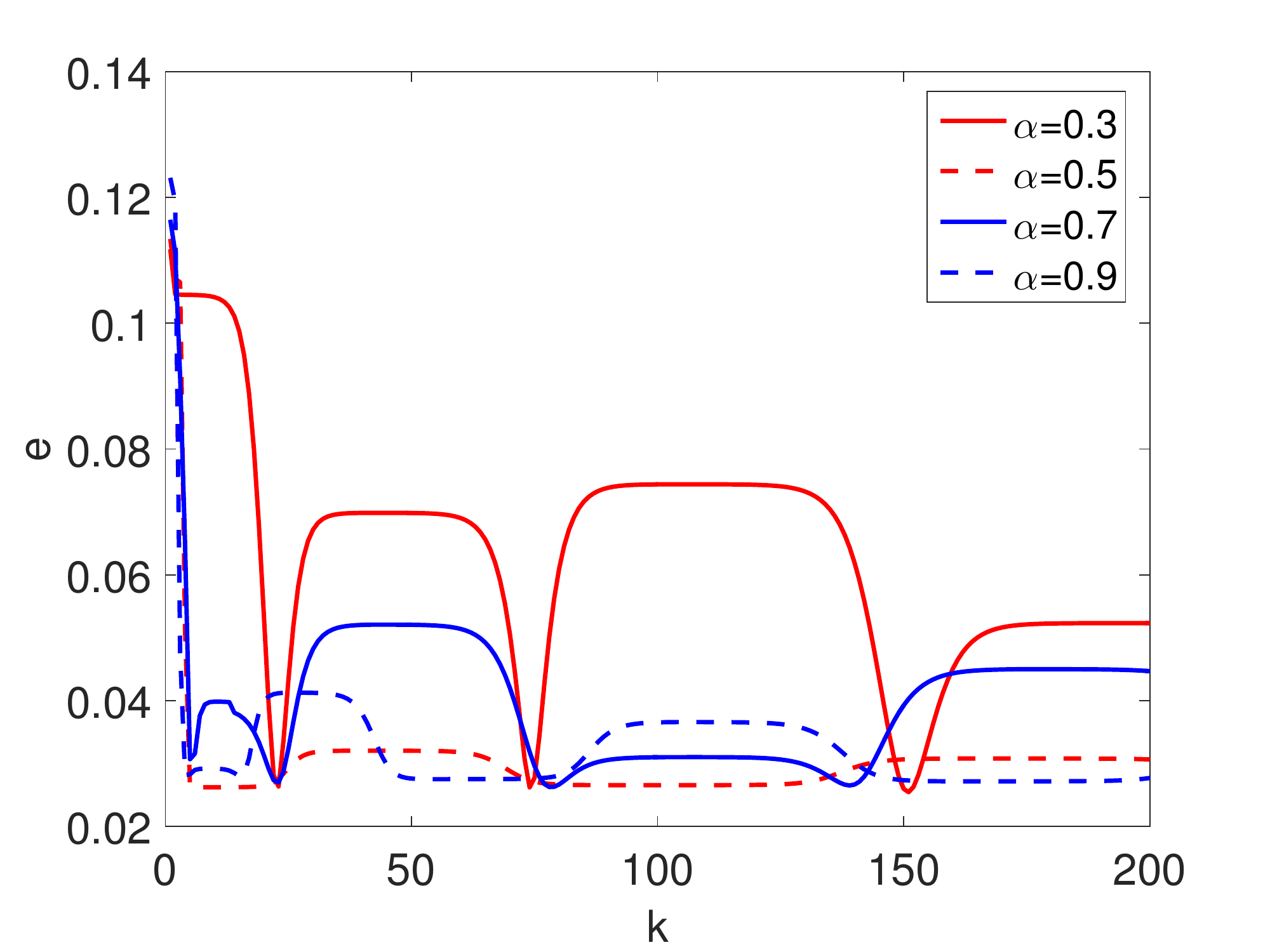} & \includegraphics[width=.33\textwidth]{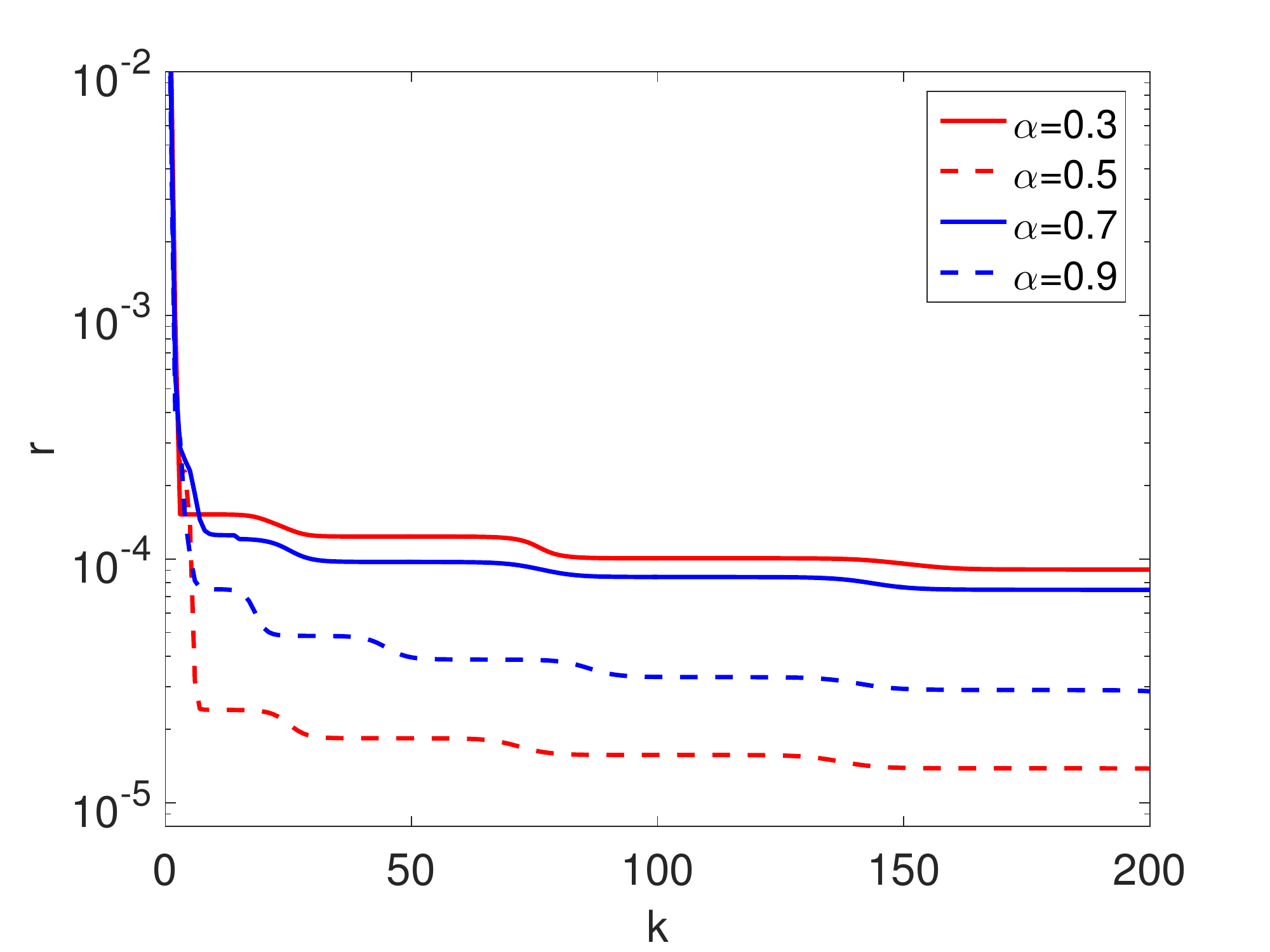}\\
 (a) reconstructions & (b) error & (c) residual
\end{tabular}
\caption{The reconstructions of the potential $q$ and the convergence behavior of the CG method in terms of the error $e$ and residual $r$ for cases (i) (top) and (ii) (bottom).\label{fig:err:q}}
\end{figure}

The reconstructions of the potential $q$ by the conjugate gradient (CG) method, based on the reduced data $\bar h(t)$,
are shown in Fig. \ref{fig:err:q} (with exact order) and Table \ref{tab:dal}. The maximum number of CG
iterations is fixed at $200$, and it is stopped so that the error is smallest possible. Throughout, for a
reconstruction $\hat q$, we measure the residual $r(\hat q)$ and the $L^2$ error $e(\hat q)$, defined
respectively by
\begin{equation*}
r(\hat q) = \|F(\hat q)-\bar h\|_{L^2(T_1,T)} \quad \mbox{and}\quad e(\hat q) = \|\hat q-q^\dag\|_{L^2(\Omega)},
\end{equation*}
where $q^\dag$ denotes the exact potential. The accuracy of the reconstructions actually does not depend
on very much on the order $\alpha$, and all the reconstructions represent a reasonable but not perfect
approximation to the true potential $q^\dag$. This observation is
consistent with prior numerical results for similar problems \cite{SunWei:2017,RundellYamamoto:2018},
and might be attributed to severe ill-conditioning of the inverse problem. The CG method
can steadily decreases the value of the objective (i.e., the residual $r$), with the first few
steps converging fairly rapidly and then slowing down considerably. Nonetheless, the error
$e$ trajectory exhibits an unusual oscillating pattern during the iteration: the error $e$ first decreases, and then
increases and then further decreases again, and there is also a flat region for which the error $e$
stays nearly constant. This behavior differs drastically from the typical steady error convergence
observed for other inverse problems, e.g., inverse source problems  \cite{JinKianZhou:2021}.
The precise mechanism of the behavior remains elusive. It is worth noting that all
these changes occur after the residual $r$ reaches a relatively small magnitude (and flat region),
indicating a potential numerical ``identifiability'' issue, despite the uniqueness in theorem
\ref{thm:main}. This also indicates that in the presence of data noise, the magnitude of the
noise has to be very small so that not to wash away these tiny transitions in order to have a fair recovery.

In the current context, the order $\alpha$ is numerically recovered, which incurs inevitable errors.
This error can potentially impact the subsequent inversion of the potential $q$. To examine the influence, we
perturb the order $\alpha$ in the optimization problem \eqref{eqn:Tikh} by
$\delta\alpha$, and repeat the numerical experiments with $\alpha+\delta\alpha$.
The results are summarized in Table \ref{tab:dal}, where $k^*$ denotes the iteration index
at which the error is smallest, and $e^*$ and $r^*$ denote the corresponding error and residual. The
presence of perturbation $\delta\alpha$ does not affect very much the attainable accuracy,
although the error $e^*$ increases steadily with the perturbation $\delta\alpha$; and generally it
takes fewer iterations to reach the optimal accuracy. This observation is largely valid for both cases
with all fractional orders under consideration.

\begin{table}[hbt!]
\centering
\caption{Numerical results for recovering $q$ (with different $\delta\alpha$).\label{tab:dal}}
\begin{threeparttable}
   \subfloat[case (i)]{
 \begin{tabular}{cccccccccc}
  \toprule
  \multicolumn{1}{c}{} & \multicolumn{3}{c}{$\delta\alpha=0$} &\multicolumn{3}{c}{$\delta\alpha=0.001$} & \multicolumn{3}{c}{$\delta\alpha=0.005$}\\
  \cmidrule(lr){2-4} \cmidrule(lr){5-7} \cmidrule(lr){8-10}
  $T_0$ & $e^*$ & $k^*$ & $r^*$ & $e^*$ & $k^*$ & $r^*$ & $e^*$ & $k^*$ & $r^*$ \\
  \midrule
  0.3 & 1.73e-2 &  98 &  2.44e-5  &  3.33e-2 &  5  &  2.35e-4 &  5.40e-2 &  2  &  1.04e-3\\
  0.5 & 1.78e-2 &  104&  2.64e-5  &  2.24e-2 &  4  &  2.63e-4 &  5.21e-2 &  3  &  1.04e-3\\
  0.7 & 1.86e-2 &  6  &  5.87e-5  &  3.68e-2 &  41 &  8.60e-5 &  2.00e-2 &  5  &  1.04e-3\\
  0.9 & 2.04e-2 &  7  &  2.62e-5  &  2.26e-2 &  41 &  1.59e-4 &  2.36e-2 &  11 &  1.04e-3\\
  \bottomrule
  \end{tabular}}

  \subfloat[case (ii)]{
  \begin{tabular}{cccccccccc}
  \toprule
  \multicolumn{1}{c}{} & \multicolumn{3}{c}{$\delta\alpha=0$} &\multicolumn{3}{c}{$\delta\alpha=0.001$} & \multicolumn{3}{c}{$\delta\alpha=0.005$}\\
  \cmidrule(lr){2-4} \cmidrule(lr){5-7} \cmidrule(lr){8-10}
  $T_0$ & $e^*$ & $k^*$ & $r^*$ & $e^*$ & $k^*$ & $r^*$ & $e^*$ & $k^*$ & $r^*$ \\
  \midrule
  0.3 & 2.54e-2 &  151 &  9.52e-5 &  2.60e-2 &  8  &  3.09e-4 &  4.66e-2 &  5 &  9.99e-4\\
  0.5 & 2.62e-2 &  8   &  2.40e-5 &  3.21e-2 &  39 &  1.17e-4 &  1.10e-1 &  2 &  1.34e-3\\
  0.7 & 2.62e-2 &  78  &  8.91e-5 &  2.73e-2 &  3  &  4.67e-4 &  2.78e-2 &  3 &  1.23e-3\\
  0.9 & 2.71e-2 &  163 &  2.90e-5 &  2.92e-2 &  56 &  1.77e-4 &  3.76e-2 &  73&  8.15e-4\\
  \bottomrule
  \end{tabular}}
  \end{threeparttable}
\end{table}

\begin{figure}[hbt!]\setlength{\tabcolsep}{0pt}
\begin{tabular}{ccc}
 \includegraphics[width=.33\textwidth]{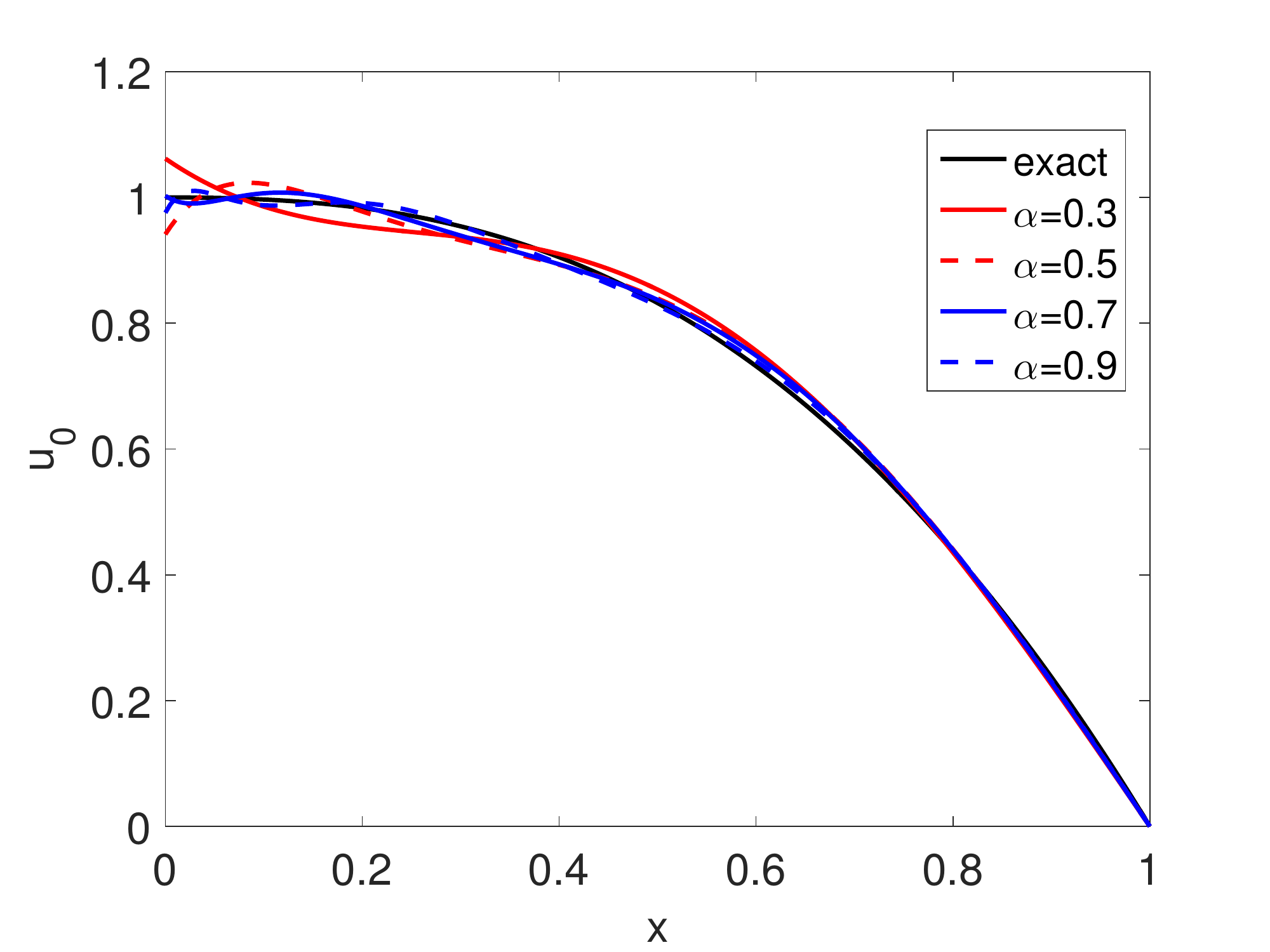} & \includegraphics[width=.33\textwidth]{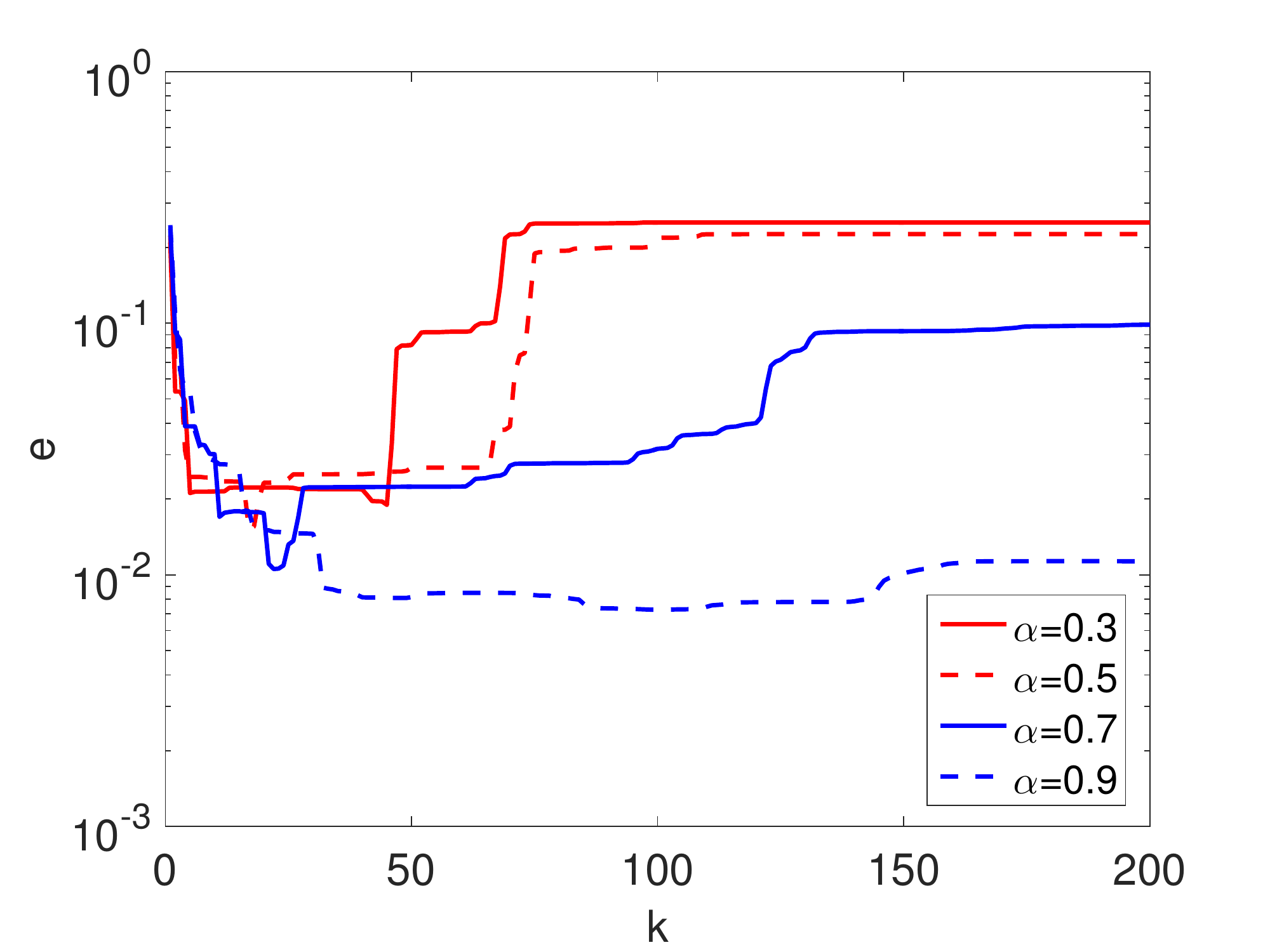} & \includegraphics[width=.33\textwidth]{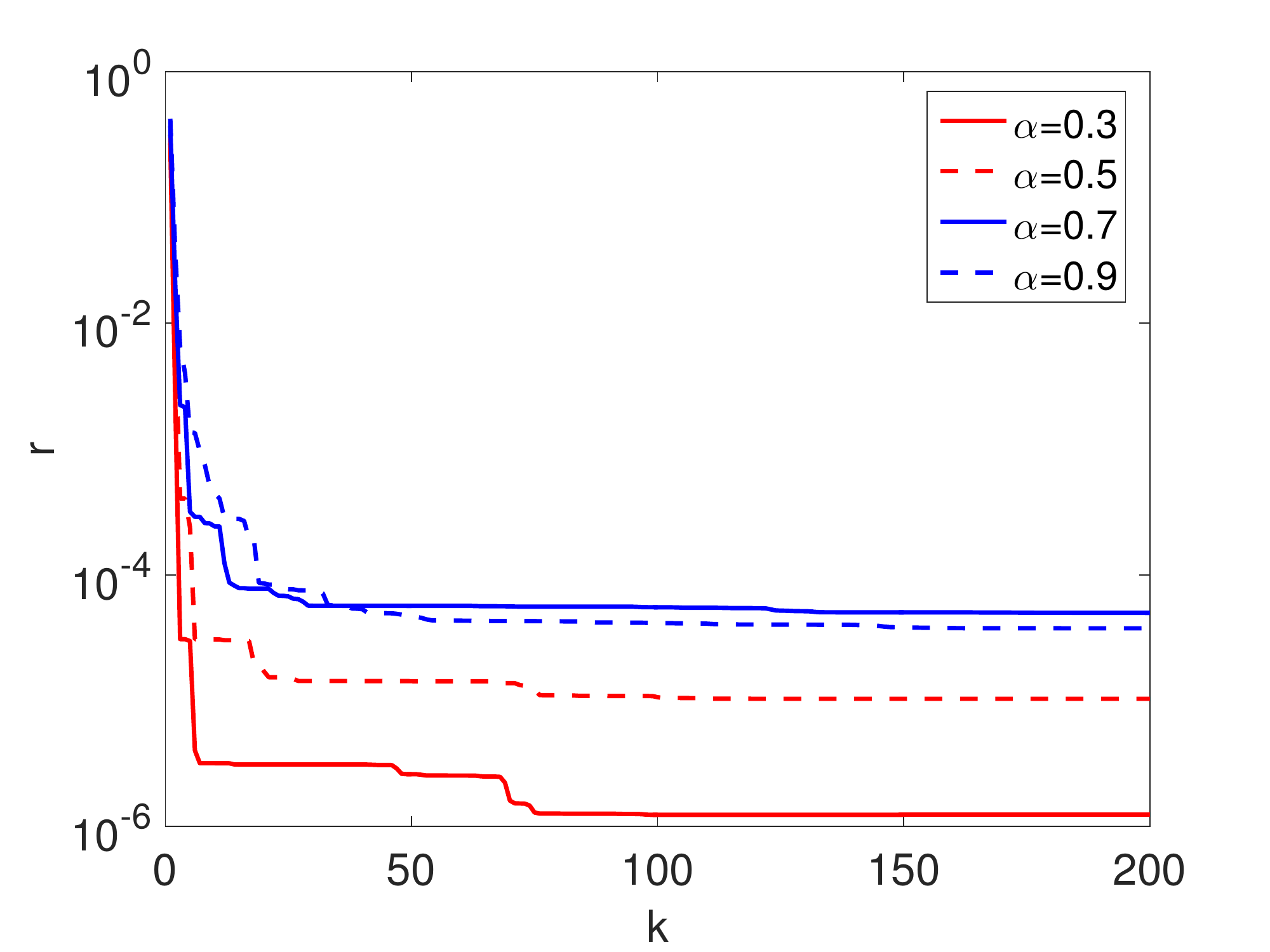}\\
  \includegraphics[width=.33\textwidth]{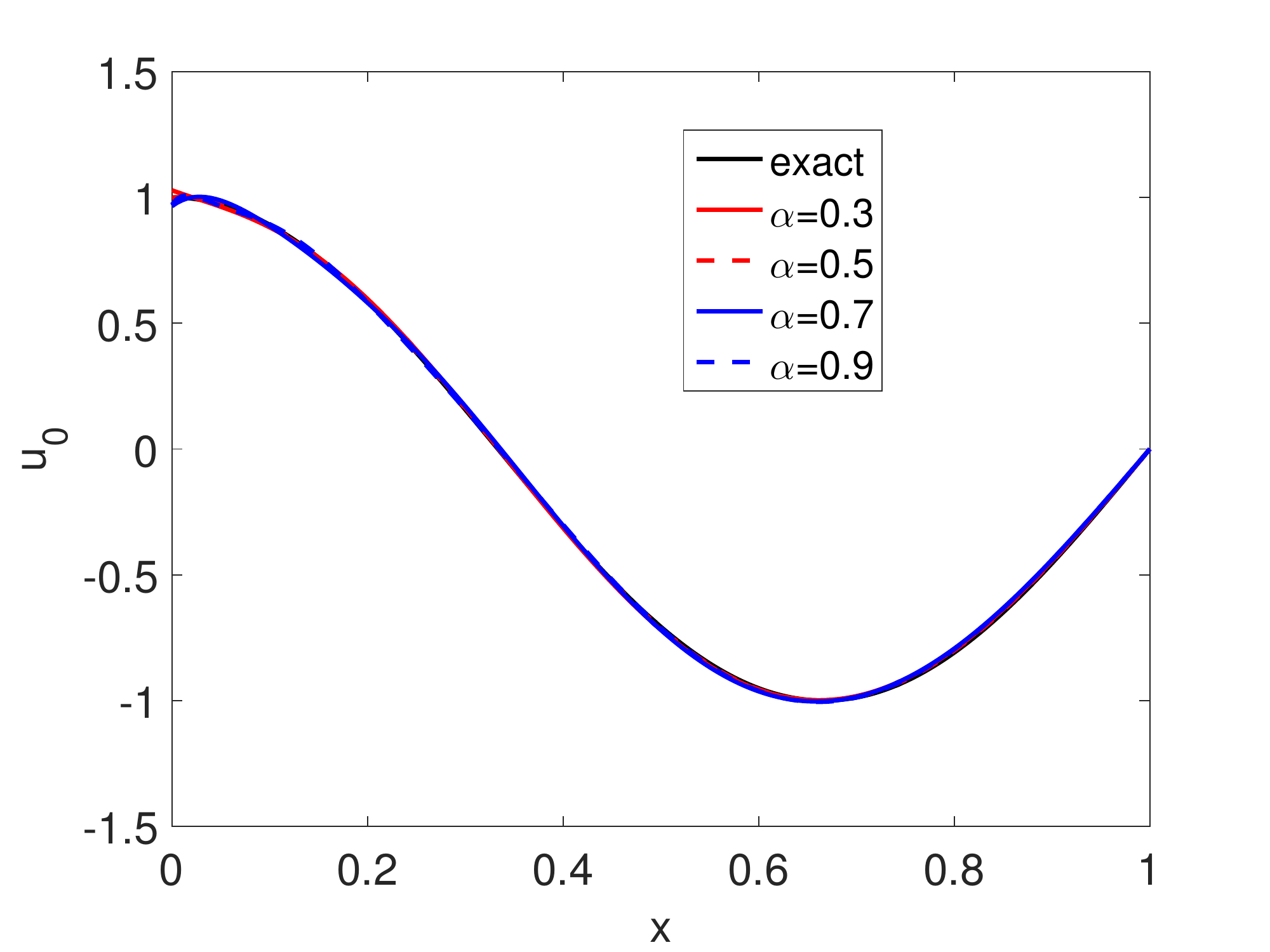} & \includegraphics[width=.33\textwidth]{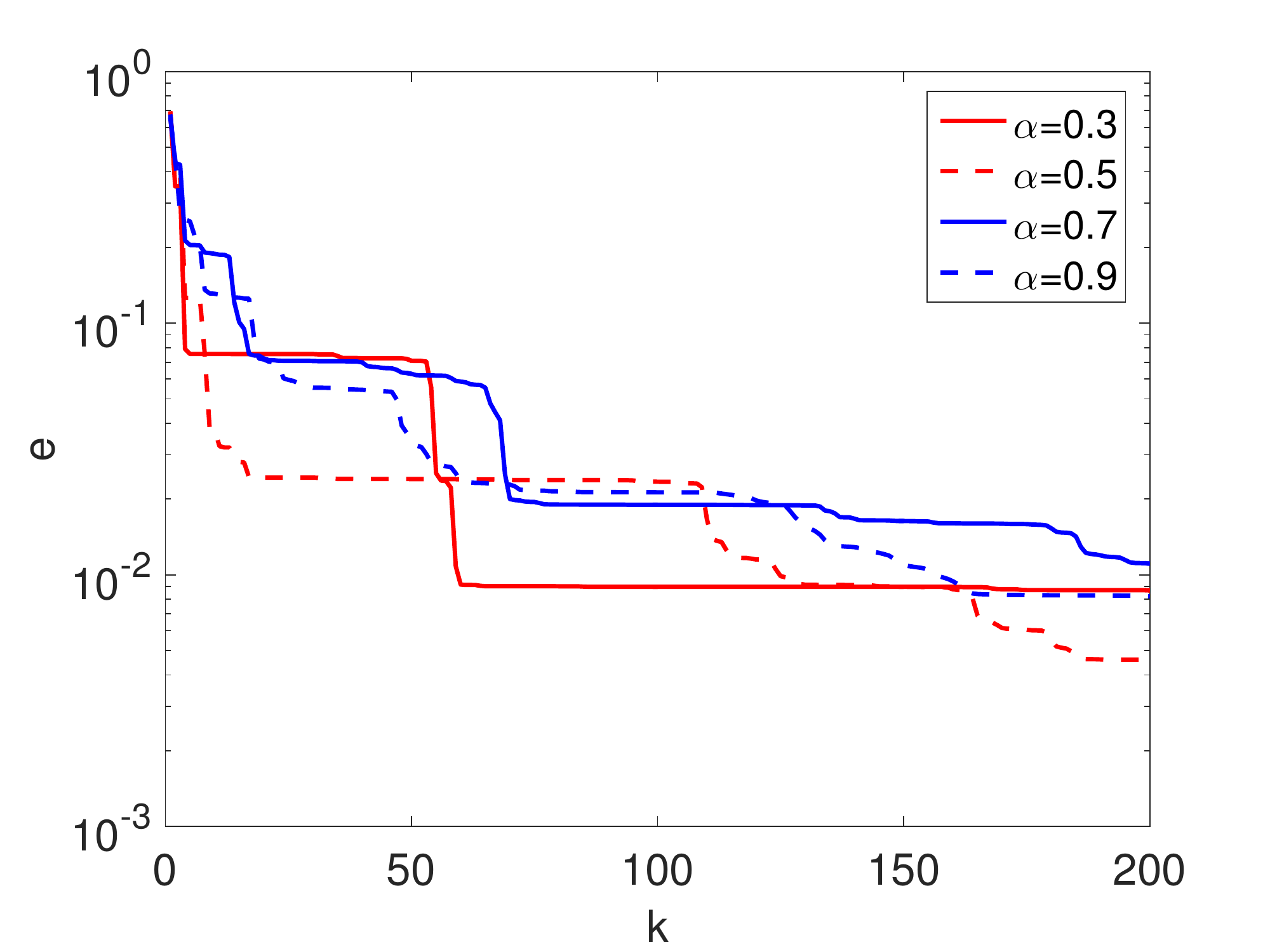} & \includegraphics[width=.33\textwidth]{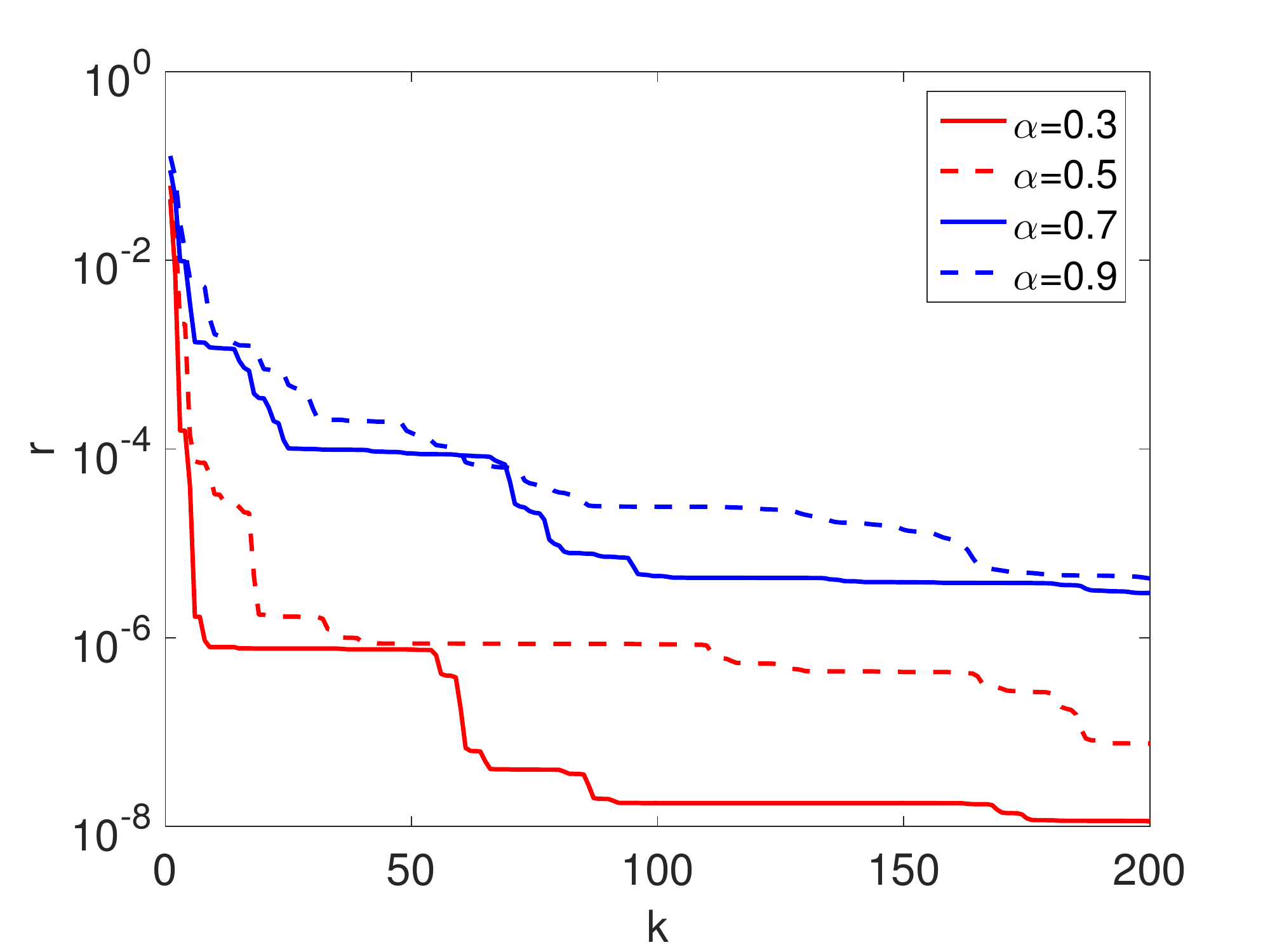}\\
 (a) reconstructions & (b) error & (c) residual
\end{tabular}
\caption{The reconstructions of the initial data $u_0$ and the convergence behavior of the CG method in terms of the error $e$ and residual $r$ for cases (i) (top) and (ii) (bottom).\label{fig:err:u0}}
\end{figure}

Last, we examine the recovery of the initial data $u_0$, using the recovered potential $q$
by the CG method (terminated after 200 iterations).
The related numerical results are summarized in Table \ref{tab:u0} and Fig. \ref{fig:err:u0},
where we have assumed that the order $\alpha$ has been estimated reliably. Interestingly,
despite the inaccuracy of the recovered potential $q$, the initial data $u_0$ can still be recovered
with a good accuracy, for all fractional orders. Further, the convergence behavior of the
CG method agrees well with that for other inverse problems (but contrasts sharply with that for
potential recovery): the method decreases the residual $e$ steadily, and the reconstruction error
$e$ exhibits a typical semiconvergence phenomenon, i.e., the iterates first converge and
then diverge, due to the inherent ill-posed nature of the inverse problem.

\begin{table}[hbt!]
\centering
\caption{Numerical results for recovering $u_0$ (inversion with recovered $q$).\label{tab:u0}}
\begin{threeparttable}
 \begin{tabular}{cccccccccc}
  \toprule
  \multicolumn{1}{c}{} & \multicolumn{3}{c}{case (i)} &\multicolumn{3}{c}{case (ii)}\\
  \cmidrule(lr){2-4} \cmidrule(lr){5-7}
  $T_0$ & $e^*$ & $k^*$ & $r^*$ & $e^*$ & $k^*$ & $r^*$  \\
  \midrule
  0.3 & 1.89e-2 &  45 & 3.07e-6 & 8.67e-3 & 200 & 1.13e-8\\
  0.5 & 1.51e-2 &  17 & 2.97e-5 & 4.59e-3 & 200 & 7.55e-8\\
  0.7 & 1.05e-2 &  22 & 7.17e-5 & 1.10e-2 & 200 & 2.97e-6\\
  0.9 & 7.26e-3 &  99 & 4.14e-5 & 8.20e-3 & 200 & 4.24e-6\\
  \bottomrule
  \end{tabular}
  \end{threeparttable}
\end{table}

In sum, the numerical experiments demonstrate the following empirical observations: (1) The order
$\alpha$ can be recovered from the observation data by a least-squares procedure; (2) the rational
function approach represents a simple method for analytically continuate the data; (3) The CG method
can produce fair approximations to the potential $q$, even under small perturbations of the order
$\alpha$, partly confirming the feasibility of the recovery, but the convergence behavior of the
algorithm exhibits an unusual oscillating feature that remains to be further examined; (4) the CG
method can produce good approximations of the initial data $u_0$, based on the recovered potential $q$.
In particular, the experiments show that the simultaneous recovery of the order, potential, and initial
data (or source) is indeed feasible, provided that accurate lateral Cauchy data is available, thereby
corroborating the uniqueness results in section \ref{sec:unique}.

\section{Concluding remarks}
In this work, we have studied an inverse problem of simultaneously recovering the fractional order and the
space-dependent potential in a one-dimensional subdiffusion model from the observation data at
the end point, when the initial data / source is not fully known. We proved that both
order and potential can be uniquely determined, if the Neumann boundary condition
satisfies a mild condition. Further, one of the space-dependent source or initial condition can
be uniquely determined, if the other is known. The analysis lends itself to an effective
two-stage reconstruction algorithm. Numerical results also show the feasibility of the recovery.

{There are many related theoretical and numerical issues awaiting further research.
First, it is of interest to extend the results to the case of a time-dependent potential. One
obstacle in the extension is that the time-dependence of the potential precludes a direct
application of the separation of variable technique, an important tool in the current analysis.
Second, it is natural to analyze more complex subdiffusion models, e.g., multi-term and
variable orders (e.g., $\alpha(t)$, $\alpha(x)$ or $\alpha(x,t)$). We believe that the results
remain largely valid for the multi-term case. However, for variable-order models, the solution
theory is still far from complete, and substantially new analytical tools are needed. Third, the
extension to the multi-dimensional case is very challenging, and requires more data for a unique
determination, e.g., restricted Neumann-to-Dirichlet map \cite{CanutoKavian:2001} or one specially
designed excitation \cite{KianLiLiu:2020}. Fourth and last, the design and analysis of relevant
reconstruction algorithms can depart enormously from the more traditional (penalized) least-squares approach.
The latter might not be directly applicable, due to the presence of unknown problem data
(and thus the very forward model in the least-squares formulation is also unknown).}

\appendix
\section{The computation of the gradients $J'(q)$ and $J'(u_0)$}
To apply the conjugate gradient method, one has to compute the gradient. This can be done
efficiently using the adjoint technique. Below we give relevant computation details
for completeness. We have the following representations of the gradients $J'(q)$ and $J'(u_0)$.
Note that the adjoint problem for $v$ and $w$ satisfies a nonlocal terminal condition.
The notation $_tI_T^{1-\alpha} v(t)$ and $_t\kern-.5em^R\kern-.2em\partial_T^\alpha v$
are defined by \cite{Jin:2021}
\begin{equation*}
_tI_T^{1-\alpha} v(t) = \frac{1}{\Gamma(1-\alpha)}\int_t^T(s-t)^{-\alpha}v(s)\d s\quad
\mbox{and}\quad _t\kern-.5em^R\kern-.2em\partial_T^\alpha v(t) = -\frac{1}{\Gamma(1-\alpha)}\frac{\d}{\d t}\int_t^T (s-t)^{-\alpha}u(s)\d s.
\end{equation*}
\begin{proposition}
The gradients $J'(q)$ and $J'(u_0)$ are respectively given by
\begin{equation*}
  J'(q) = -\int_{0}^T u(q)v(q)\d t\quad \mbox{and}\quad
  J'(u_0) =- ({_tI_{T_1}^{1-\alpha}}w)(0) = -\frac{1}{\Gamma(1-\alpha)}\int_0^{T_1}t^{-\alpha}w(t)\d t,
\end{equation*}
with $v\equiv v(q)$ and $w$ solving respectively
\begin{equation*}
  \left\{\begin{aligned}
    _t\kern-.5em^R\kern-.2em\partial_T^\alpha v - \mathcal{A}v &= 0,\quad \mbox{in }\Omega\times [0,T),\\
    {_tI_T^{1-\alpha}}v(x,T) & = 0,\quad \mbox{in }\Omega,\\
    \partial_{x} v(0,\ell,t) & = F(q)-\bar h, \quad \mbox{in }[0,T),\\
    v(1,t) & = 0, \quad \mbox{in }[0,T),
  \end{aligned}  \right.
  \quad \mbox{and}\quad
  \left\{\begin{aligned}
    _t\kern-.5em^R\kern-.2em\partial_{T_1}^\alpha w - \mathcal{A} w &= 0,\quad \mbox{in }\Omega\times [0,T_1),\\
    {_tI_{T_1}^{1-\alpha}}w(x,T_1) & = 0,\quad \mbox{in }\Omega,\\
    \partial_{x} w(0,\ell,t) & = F(u_0)- h, \quad \mbox{in }[0,T_1),\\
    w(1,t) & = 0, \quad \mbox{in }[0,T_1].
  \end{aligned}  \right.
\end{equation*}
\end{proposition}
\begin{proof}
Let $X=\{v\in L^2(0,T;H^1(\Omega)): v(1,t)=0, t\in(0,T)\}$.
The directional derivative $J'(q)[\delta q]$ with $\delta q\in L^2(\Omega)$ is given by
$J'(q)[\delta q] = (u_q'(q)[\delta q],F(q)-\bar h)_{L^2(0,T)}$,
where $u^\delta =u_q'(q)[\delta q]$ satisfies $u^\delta(0)=0$ and
\begin{equation}\label{eqn:fde-lin}
  \int_0^T\int_\Omega (\phi\,\partial_t^\alpha u^\delta+ a\nabla u^\delta\cdot\nabla \phi +qu^\delta\phi) \d x\d t = -\int_0^T\int_\Omega \delta qu(q)\phi\d x \d t,\quad \forall \phi \in X.
\end{equation}
Meanwhile, the weak formulation for the adjoint solution $v$ is given by
\begin{equation}\label{eqn:fde-adj-weak}
  \int_0^T\int_\Omega (\phi\, {_t\kern-.5em^R\kern-.2em\partial_T^\alpha}v +a \nabla v\cdot\nabla \phi+qv\phi)\d x\d t = \int_0^T(F(q)-\bar h)\phi(0,t)\d t,\quad \forall \phi\in X.
\end{equation}
Then taking $\phi=v \in X$ in \eqref{eqn:fde-lin} and $\phi=u^\delta \in X$ in \eqref{eqn:fde-adj-weak}, using
the following integration by parts formula (see, e.g., \cite[p. 76, Lemma 2.7]{KilbasSrivastavaTrujillo:2006} or \cite[Lemma 2.6]{Jin:2021})
\begin{equation}\label{eqn:int-by-part}
  \int_0^T \psi\, {\partial_t^\alpha \phi} \d t = (\phi\, {_tI_T^{1-\alpha}}\psi)|_{t=0}^T + \int_0^T \phi\,\, { _t\kern-.5em^R\kern-.2em\partial_T^\alpha \psi}\d t,
\end{equation}
and last subtracting the two identities give
\begin{equation*}
 - \int_0^T\int_\Omega \delta qu(q)v\d x \d t = \int_0^T(F(q)-\bar h)u^\delta(0,t)\d t,
\end{equation*}
from which we deduce the expression of $J'(q)$.
Similarly, $J'(u_0)[\delta u_0] = (u_{u_0}'(u_0)[\delta u_0],F(u_0)- h)_{L^2(0,T_1)}$,
where $u^\delta =u_{u_0}'(u_0)[\delta u_0]$ (slightly abused notation) satisfies $u^\delta(0)=\delta u_0$ and
\begin{equation}\label{eqn:fde-lin0}
  \int_0^{T_1}\int_\Omega (\phi\,\partial_t^\alpha u^\delta+ a\nabla u^\delta\cdot\nabla \phi +qu^\delta\phi) \d x\d t = 0,\quad \forall \phi \in X.
\end{equation}
Meanwhile, the space-time weak formulation for the adjoint solution $w$ is given by
\begin{equation}\label{eqn:fde-adj-weak0}
  \int_0^{T_1}\int_\Omega (\phi\, {_t\kern-.5em^R\kern-.2em\partial_T^\alpha}w + a\nabla w\cdot\nabla \phi+qw\phi)\d x\d t = \int_0^{T_1}(F(u_0)- h)\phi(0,t)\d t,\quad \forall \phi \in X.
\end{equation}
Then taking $\phi=w$ in \eqref{eqn:fde-lin0} and $\phi=u^\delta$ in \eqref{eqn:fde-adj-weak0},
applying the integration by parts formula \eqref{eqn:int-by-part}, and subtracting the two identities give
\begin{equation*}
 - \int_0^{T_1}\int_\Omega \delta u_0 {_tI_{T_1}^{1-\alpha}}w(t)\d x \d t = \int_0^{T_1}(F(q)- h)u^\delta(0,t)\d t.
\end{equation*}
This gives the expression of $J'(u_0)$.
\end{proof}

\begin{remark}
One can also derive the regularity of the gradients.
For example, with $g\in L^2(0,T)$, we have $u\in L^2(0,T;D(A^s))$, for any
$s\in (\frac12,\frac34)$, cf. the proof of proposition \ref{prop:solrep}. Similarly
for $\bar h\in L^2(0,T)$, there holds the adjoint $v\in L^2(0,T;D(A^s))$. This and
algebraic property of the space $D(A^s)$ imply $uv\in L^1(0,T;D(A^s))$, and thus
$J'(q)\in D(A^s)$.
\end{remark}

\bibliographystyle{abbrv}
\bibliography{frac}
\end{document}